\documentclass[10pt]{amsart}

\usepackage{amsmath}
\usepackage{amsthm}
\usepackage{amssymb}
\usepackage{amsfonts}  
\usepackage[all]{xy}

\newcommand{\cala}{{\mathcal{A}}}

\newcommand{\calc}{{\mathcal{C}}}
\newcommand{\cald}{{\mathcal{D}}}
\newcommand{\cale}{{\mathcal{E}}}
\newcommand{\calf}{{\mathcal{F}}}

\newcommand{\cals}{{\mathcal{S}}}
\newcommand{\calt}{{\mathcal{T}}}

\newcommand{\calw}{{\mathcal{W}}}

\newcommand{\IB}{{\mathbb B}}

\newcommand{\IH}{{\mathbb H}}

\newcommand{\IK}{{\mathbb K}}

\newcommand{\IN}{{\mathbb N}}

\newcommand{\IR}{{\mathbb R}} 

\newcommand{\IT}{{\mathbb T}}

\newcommand{\IZ}{{\mathbb Z}} 

\newcommand{\Jt}{\mathbf{t}}
\newcommand{\Jd}{\mathbf{d}}

\newcommand{\All}{{\mathit{All}}} 
\newcommand{\Cyc}{{\mathit{Cyc}}} 
\newcommand{\VCyc}{{\mathit{VCyc}}} 
\newcommand{\mor}{{\rm mor}} 
\newcommand{\supp}{{\rm supp}} 
\newcommand{\sma}{{\wedge}} 
\newcommand{\einsu}{{\left[ 1, \infty \right)}}

\newcommand{\inc}{{\rm inc}}

\newcommand{\id}{{\rm id}}
\newcommand{\tr}{{\rm tr}}

\newcommand{\pt}{{\rm pt}}

\newcommand{\cone}{{\rm cone}}
\newcommand{\cyl}{{\rm cyl}}
\newcommand{\ch}{{\rm ch}}      
\newcommand{\NK}{{\rm NK}}      
\newcommand{\sh}{{\mathit{sh}}}
\newcommand{\step}{{\rm step}}
\newcommand{\hocolim}{\operatorname{hocolim}}

\newcommand{\fr}{{\mathit{fr}}}
\newcommand{\hf}{{\mathit{hf}}}

\newcommand{\sing}{{\mathit{sing}}}

\newcommand{\Sw}{{\mathit{Sw}}}

\newcommand{\e}{\epsilon}

\newcommand{\x}{\times}
\newcommand{\dd}{\partial}

\newcommand{\asy}{\mathit{asy}}

\theoremstyle{plain}
\newtheorem{theorem}{Theorem}[section]
\newtheorem{lemma}[theorem]{Lemma}
\newtheorem{corollary}[theorem]{Corollary}
\newtheorem{proposition}[theorem]{Proposition}
\newtheorem{conjecture}[theorem]{Conjecture}

\theoremstyle{definition}
\newtheorem{definition}[theorem]{Definition}
\newtheorem{example}[theorem]{Example}

\newtheorem{warning}[theorem]{Warning}

\theoremstyle{remark}
\newtheorem{remark}[theorem]{Remark}

\newtheorem{assumption}[theorem]{Assumption}

\makeatletter\let\c@equation=\c@theorem\makeatother

\begin{document}

\typeout{---------------------- iso2 -----------------------------}

\begin{titlepage} 
\title{On the Farrell-Jones Conjecture for higher algebraic K-Theory} 
\author{Arthur Bartels, Holger Reich}
\address{Fachbereich Mathematik\\ Universit{\"a}t M{\"u}nster \\ Einsteinstr.~62\\48149~M{\"u}nster\\Germany}
\email{bartelsa@math.uni-muenster.de}
\email{reichh@math.uni-muenster.de}
\date{\today}
\begin{abstract}
We prove the Farrell-Jones Conjecture about the algebraic $K$-theory of a group ring $R \Gamma$ in 
the case where the group $\Gamma$ is the fundamental group of a closed Riemannian manifold with strictly negative sectional
curvature. The coefficient ring $R$ is an arbitrary associative ring with unit and the result applies to all dimensions.
\end{abstract}
\end{titlepage}

\maketitle

\typeout{--------------------- Inhalt ----------------------------}
\tableofcontents


\typeout{---------------------einfuehrung-------------------------}

\section{Introduction}

Conjecturally the algebraic $K$-theory groups $K_n( \IZ \Gamma )$, $n \in \IZ$, of the integral group ring $\IZ \Gamma$ of every  torsion free
group $\Gamma$ can be expressed in terms of the group homology of $\Gamma$ and the algebraic $K$-theory of the integers.
More precisely there is the following conjecture, compare \cite[Section~VI]{Hsiang(icm83)}.
\begin{conjecture} \label{con:Hsiang-conjecture}
For a torsion free group $\Gamma$ the so called assembly map \cite{Loday(assembly)}
\[
A \colon H_n ( B\Gamma ; \IK^{-\infty} ( \IZ )) \to K_n ( \IZ \Gamma )
\]
is an isomorphism for all $n \in \IZ$.
\end{conjecture}
Here $B\Gamma$ is the classifying space of the group $\Gamma$ and we denote 
by $\IK^{-\infty} ( R )$ the non-connective algebraic $K$-theory spectrum 
of the ring $R$. The homotopy groups of this  spectrum are denoted $K_n ( R )$ and coincide with 
Quillen`s algebraic $K$-groups of $R$ in positive dimensions \cite{Quillen(higher-K)} and with the negative 
$K$-groups of Bass \cite{Bass(book)} in negative dimensions. 
The homotopy groups of the spectrum $X_+ \sma \IK^{-\infty} ( \IZ )$ are denoted $H_n ( X ; \IK^{-\infty} ( \IZ ))$.
They yield  a generalized homology theory and in particular standard computational tools like the Atiyah-Hirzebruch spectral sequence
apply to the left hand side of the assembly map above.

As a corollary of the main result of this paper we prove Conjecture~\ref{con:Hsiang-conjecture}
in the case where $\Gamma$ is the fundamental group of a closed Riemannian
manifold with strictly negative sectional curvature.
In fact our result is more general and applies to group rings $R\Gamma$, where $R$ is a completely arbitrary
coefficient ring.

Note that if 
one replaces in Conjecture~\ref{con:Hsiang-conjecture} 
the coefficient ring  $\IZ$ by an arbitrary coefficient ring $R$ 
the corresponding conjecture would be false already in the simplest non-trivial case: if $\Gamma=C$ is the infinite cyclic group
the Bass-Heller-Swan formula 
\cite{Bass-Heller-Swan(1964)}, 
\cite[p.236]{Grayson(higher-K-II)} for $K_n( RC ) = K_n ( R [t^{\pm 1 }] )$ yields that 
\[
K_n ( R C ) \cong K_{n-1} ( R ) \oplus K_n ( R ) \oplus \NK_n ( R ) \oplus \NK_n ( R ),
\]
where $\NK_n ( R )$ is defined as the cokernel of the split inclusion $K_n ( R ) \to K_n ( R [t] )$ and does not vanish
in general. But since $S^1$ is a model for $BC$ one obtains on the left hand side of the assembly map only
\[
H_n ( BC ; \IK^{-\infty}( R ) ) \cong K_n ( R) \oplus K_{n-1} (R).
\]

In some sense this is all that goes wrong.
Combining our main result Theorem~\ref{main-theorem} with Proposition~1.8 in \cite{Bartels-Farrell-Jones-Reich(topology)}
we obtain the following generalization of the Bass-Heller-Swan formula.
\begin{corollary} \label{cor:generalized-BHS}
Let $\Gamma$ be a fundamental group of a closed  Riemannian manifold with strictly negative sectional curvature, then
for every associative ring with unit $R$ we have
\[
K_n ( R \Gamma ) \cong H_n ( B \Gamma ; \IK^{-\infty} ( R ) ) \oplus \bigoplus_{I} \NK_n ( R ) \oplus \NK_n ( R )
\]
for all $n \in \IZ$.
Here the sum on the right is indexed over the set $I$ of conjugacy classes of maximal cyclic subgroups of $\Gamma$.
\end{corollary}
Recall that a ring $R$ is called (right) regular if it is (right) Noetherian and every finitely generated (right) 
$R$-module admits a finite dimensional projective resolution. Principal ideal domains are examples of regular rings.
It is known \cite[Chapter~XII]{Bass(book)}, \cite[p.122]{Quillen(higher-K)}  
that for a regular ring $\NK_n(R)=0$ for all $n \in \IZ$. Hence for regular coefficient rings 
the expression in Corollary~\ref{cor:generalized-BHS}
simplifies and proves the more general version of Conjecture~\ref{con:Hsiang-conjecture}
where the coefficient ring $\IZ$ is replaced by an arbitrary regular coefficient ring $R$.

We proceed to describe the Farrell-Jones Conjecture about algebraic $K$-theory of group rings \cite{Farrell-Jones(isomorphism)} which is 
the correct conceptional framework for these kinds of results and which applies also to groups which contain torsion. 
Here is some more notation.

A set of subgroups of a given group $\Gamma$ is called a family of subgroups if 
it is closed under conjugation with elements from $\Gamma$ and closed under taking subgroups.
We denote by 
\[
\{ 1 \} , \quad \Cyc, \quad \VCyc \quad \mbox{ and } \quad \All
\]
the families which consist of the trivial subgroup, all cyclic subgroups, all virtually cyclic subgroups, respectively 
all subgroups of $\Gamma$.
Recall that a group is called virtually cyclic if it contains a cyclic subgroup of finite index. 

For every family $\calf$ of subgroups of $\Gamma$ there exists a classifying space for the family $\calf$ denoted
$E \Gamma ( \calf )$, 
compare~\cite{tomDieck(orbittypen)}, 
\cite[I.6]{tomDieck(transformation-groups)}, 
and \cite[Appendix]{Farrell-Jones(isomorphism)}. 
It is characterized by the universal property that for every $\Gamma$-CW-complex $X$ 
whose isotropy groups are all in the family $\calf$ there exists an equivariant continuous map $X \to E \Gamma ( \calf )$
which is unique up to equivariant homotopy. A $\Gamma$-CW-complex $E$ is a model for the classifying space $E \Gamma ( \calf )$
if the fixpoint sets $E^H$ are contractible for $H \in \calf $ and empty otherwise. Note that the one point space $\pt$
is a model for $E \Gamma ( \All )$ and that $E \Gamma ( \{ 1 \} )$ is the universal covering of the classifying space $B \Gamma$.

In \cite{Davis-Lueck(assembly)} Davis and L{\"u}ck construct a generalized equivariant homology theory for $\Gamma$-CW-complexes
$X \mapsto H_n^{\Gamma} ( X ; \IK_R^{-\infty} )$
associated to a jazzed-up version of the non-connective algebraic $K$-theory spectrum functor. If one evaluates this 
$\Gamma$-homology theory on a homogeneous space $\Gamma / H $ one obtains
$H_n^{\Gamma} ( \Gamma / H ; \IK^{-\infty}_R ) \cong K_n ( R H )$.
Using this language the Farrell-Jones Conjecture for the algebraic $K$-theory of group rings can be formulated as follows.

\begin{conjecture}[The Farrell-Jones Conjecture for $K_{n} ( R \Gamma )$] \label{con:FJcon}
Let $\Gamma$ be a group and let $R$ be an associative ring with unit. Then the map
\[
A_{\VCyc} \colon H_n^{\Gamma} ( E \Gamma ( \VCyc ) ; \IK^{-\infty}_R ) \to H_n^{\Gamma} ( \pt  ; \IK^{-\infty}_R ) \cong K_n ( R \Gamma ).
\]
which is induced by the projection $E \Gamma ( \VCyc ) \to \pt$  is an isomorphism for all $n \in \IZ$.
\end{conjecture}

This conjecture was formulated in \cite{Farrell-Jones(isomorphism)} for $R= \IZ$ and stated in this more general form in 
\cite{Bartels-Farrell-Jones-Reich(topology)}. Above we used the language developed by Davis and L{\"u}ck in \cite{Davis-Lueck(assembly)}
to formulate the conjecture. The identification of this formulation with the original formulation
in \cite{Farrell-Jones(isomorphism)} which used \cite{Quinn(ends-II)} is carried out in \cite{Hambleton-Pedersen(identifying)}. 
For more information on this and related conjectures the reader should consult \cite{Lueck-Reich(survey)}.

Our main result proves Conjecture~\ref{con:FJcon} for the class of groups that was already mentioned above. 
\begin{theorem} \label{main-theorem}
Let $R$ be an associative ring with unit. Let $\Gamma$ be the fundamental group
of a closed Riemannian manifold with strictly negative sectional curvature.
Then for all $n \in \IZ$ the assembly map
\[
A_{\Cyc} \colon H_n^{\Gamma} ( E \Gamma ( \Cyc ) ; \IK^{-\infty}_R  ) \to K_n ( R \Gamma )
\]
induced by the projection $E \Gamma ( \Cyc ) \to \pt$ is an isomorphism.
\end{theorem}
Note that the groups considered in this theorem are torsion free and that  for a torsion free group the family 
$\VCyc$ of all virtually cyclic subgroups reduces to the family $\Cyc$ of all cyclic subgroups.
The result extends the results of \cite{Bartels-Farrell-Jones-Reich(topology)}, 
where surjectivity in low dimensions and injectivity was proven.
Results which are strongly related to the result above about the low dimensional $K$-theory of 
the integral group ring, pseudoisotopy spectra and the structure set in surgery theory were 
proven by Farrell and Jones in \cite{Farrell-Jones(dynamics-I)}, 
\cite{Farrell-Jones(dynamics-II)}, 
\cite{Farrell-Jones(borel)}
and 
\cite{Farrell-Jones(pseudo-non-positive)}. 
Apart from \cite{Waldhausen(generalized-free-I)} the result above seems to be the first integral result of this 
type which applies to the higher algebraic $K$-theory of group rings.

From the fact that we do not make any assumptions on the coefficient ring $R$ one can derive a corresponding isomorphism
statement for an assembly map for $\NK$-groups.
\begin{corollary}
Let $R$ be an associative ring with unit. Let $\Gamma$ be the fundamental group of closed Riemannian manifold with 
strictly negative sectional curvature. Then for all $n \in \IZ$ the assembly map for $\NK$-groups
\[
A_{\Cyc} \colon H_n^{\Gamma} ( E \Gamma ( \Cyc ) ; \IN \IK^{-\infty}_R  ) 
\to H_n^{\Gamma} ( \pt ; \IN \IK^{-\infty}_R ) \cong \NK_n ( R \Gamma )
\]
is an isomorphism.
\end{corollary}
\begin{proof} Since there is a splitting $\IK_{R[t]}^{-\infty} \simeq \IK_{R}^{-\infty} \vee \IN \IK_{R}^{-\infty}$ 
the isomorphism result for two of the assembly maps associated to $\IK_R^{-\infty}$, $\IK_{R[t]}^{-\infty}$, respectively 
$\IN \IK_{R}^{-\infty}$ implies it for the third, compare~\cite[Proposition~7.4]{Bartels-Farrell-Jones-Reich(topology)}.
\end{proof}
In particular one can conclude that the vanishing of $\NK_n ( R )$ for 
$n \leq N$ implies the vanishing of $\NK_n (R \Gamma )$ for $n \leq N$.
Note that even if $R$ is regular it is not at all clear if $R \Gamma$ is regular.
If one uses that $R [ \Gamma_1 \times \Gamma_2] = R[\Gamma_1][\Gamma_2]$ and 
iterates one obtains the following corollary.
\begin{corollary} \label{cor:iterate-corollary}
Suppose $\Gamma = \Gamma_1 \times \Gamma_2 \times \dots \times \Gamma_k$, where each $\Gamma_i$ is the fundamental group
of a closed Riemannian manifold with strictly negative sectional curvature. If $R$ is a regular ring, then
the assembly map
\[
A \colon H_n ( B \Gamma ; \IK^{-\infty} ( R ) ) \to K_n ( R \Gamma )
\]
is an isomorphism for all $n \in \IZ$.
\end{corollary}

The proof of Theorem~\ref{main-theorem} relies on the fact 
(see Subsection~\ref{subsec:assembly-as-forget-control} and \cite{Bartels-Farrell-Jones-Reich(topology)}) 
that the assembly map $A_{\VCyc}$ can be described as a ``forget-control map''
in the sense of controlled topology, compare \cite{Quinn(ends-II)}, \cite{Pedersen-Weibel(homology)}.
In order to  prove a surjectivity result we hence have to ``gain control''.  More precisely 
the negatively curved manifold $M$ whose fundamental group we want to treat can 
be used in order to construct a geometric model for the map $E \Gamma \to E \Gamma ( \Cyc )$ given by the universal
property. We will consider 
suitable additive categories of $R$-modules and morphisms over $E \Gamma \times [ 1 , \infty)$ where the morphisms satisfy
control conditions. The assembly map is obtained by applying $K$-theory to an inclusion of additive categories, where the larger category
differs from the smaller one by a relaxed control condition on the morphisms.

The geometric program in order to gain control stems from \cite{Farrell-Jones(dynamics-I)} and consists mainly of three steps:
\begin{enumerate}
\item[(I)]
Construct a transfer from the manifold $M$ to a suitable subbundle of its sphere bundle $SM$ in such a way that
transferred morphisms are in asymptotic position, i.e.\ they are in a good starting position for the geodesic flow.
Make sure that transferring up and projecting down again does not change the $K$-theory class.
\item[(II)]
Consider the foliation on the sphere bundle $SM$ given by the flow lines of the geodesic flow. 
Use the geodesic flow on $SM$ in order to achieve foliated control.
\item[(III)]
Prove a ``Foliated Control Theorem'' in order to improve from ``foliated control'' to ``ordinary control''.
At least do so away from the ``short'' closed geodesics.
\end{enumerate}
Note that the closed geodesics that appear in Step~(III) are in bijection  with conjugacy classes of cyclic subgroups.
Hence the family $\Cyc$ which appears in Theorem~\ref{main-theorem} shows up quite naturally in the proof.

We refer to Section~\ref{section-proof} for a more detailed outline of the proof. In the following we only discuss
why new techniques were necessary in order to treat higher algebraic 
$K$-theory along the lines of the program above.

One main difficulty was to construct a suitable transfer as required in Step~(I) of the program.
Looking at the analogous situations for $h$-cobordisms or $A$-theory, where a transfer is given by pull-back, it is in 
principle clear what the algebraic analogue in our situation should be. However the obvious naive approaches
are not ``functorial enough'' to induce a map in higher $K$-theory.
Hence one needs to come up with a suitably refined transfer which takes care of the functoriality problems
(e.g.\ work with  singular chain complexes) but at the same time does not destroy the control requirements.
In order to treat the question whether transferring up and projecting down yields the identity on $K$-theory we prove in
Proposition~\ref{mult-by-swan} a formula for the kind of transfers we construct. Transferring up and projecting down 
yields multiplication by a certain element in the Swan group. The Swan group element depends on the homology groups of the fiber 
considered as modules over the fundamental group of the base.

Another main achievement in this paper is the Foliated Control Theorem~\ref{fct} for higher algebraic
$K$-theory. Earlier foliated control theorems (see for example Theorem~1.6 in \cite{Farrell-Jones(dynamics-I)} 
or Theorem~1.1 in \cite{Bartels-Farrell-Jones-Reich(trieste)})
were formulated for individual $K$-theory elements. It is however difficult to explicitly describe elements
in higher $K$-theory groups. Hence we had to find a way to formulate and prove a foliated control theorem
in a more ``functorial'' fashion. 
We would like to emphasize that the Foliated Control Theorem~\ref{fct} does not rely on a squeezing result
or any kind of torus trick. Essentially only the existence of the long exact sequence associated to 
a Karoubi filtration \cite{Cardenas-Pedersen(karoubi)} and Eilenberg-swindles 
are used as the abstract building blocks for such an argument. 
Of course on the geometric side the existence of long-and-thin cell structures proven in \cite{Farrell-Jones(dynamics-I)} is crucial.
In particular our technique should prove analogous foliated control theorems in the context of 
algebraic $L$-theory or topological $K$-theory of $C^{\ast}$-algebras
since the corresponding tools are available in those set-ups, compare~\cite{Higson-Pedersen-Roe(controlled)}.

The reader who is familiar with the work of Farrell-Jones (in particular \cite{Farrell-Jones(pseudo-non-positive)} and 
\cite{Farrell-Jones(borel-non-positive)})
may wonder why 
we cannot weaken the assumption 
in Theorem~\ref{main-theorem} 
from strictly negative 
curvature to non-positive curvature. The reason is that
the focal transfer which  is used in \cite{Farrell-Jones(pseudo-non-positive)} and \cite{Farrell-Jones(borel-non-positive)}, 
in contrast to the asymptotic transfer 
used in this paper, is  definitely not functorial and it is hence even more 
difficult to describe a corresponding transfer for higher algebraic $K$-theory.

It also remains open whether the program can be adapted to prove cases of the Baum-Connes Conjecture~\cite{Baum-Connes-Higson(conjecture)}
or to treat algebraic $L$-theory with arbitrary coefficients. In both cases a crucial question is whether a suitable 
transfer can be constructed.
(A geometric version of an $L$-theory transfer is one of the many ingredients in \cite{Farrell-Jones(borel)}.)

\subsection{Acknowledgments}

We would like to thank Tom Farrell and Lowell Jones for many discussions on the subject.
We would also like to thank Wolfgang L{\"u}ck for help with the Swan group actions.
Our research was supported by the SFB 478 - Geometrische Strukturen in der Mathematik - M{\"u}nster.


\typeout{---------------------preliminaries -------------------------}

\section{Preliminaries}

\subsection{Conventions and notation}
In this section we briefly introduce some notation that is used throughout the proof. For more details the reader should consult
Section~2 in \cite{Bartels-Farrell-Jones-Reich(topology)}. 

\subsubsection{The functor $\IK^{-\infty}$}
We will denote by $\IK^{-\infty}$ the functor which associates to an additive category its non-connective 
$K$-theory spectrum, see \cite{Pedersen-Weibel(delooping)} or \cite{Cardenas-Pedersen(karoubi)}.
We assume that the reader is familiar with the standard properties of this functor, 
compare~\cite[Subsection~2.1]{Bartels-Farrell-Jones-Reich(topology)}.
A statement about exact functors between additive categories is true ``after applying $K$-theory'' or 
``on the level of $K$-theory'' will always mean after applying $\IK^{-\infty}$.

\subsubsection{Modules and morphisms over a space} \label{subsubsec:modules-and-morphisms}
As explained below in Subsection~\ref{subsec:assembly-as-forget-control} the assembly 
map will be described as a ``forget-control map'' between suitable additive categories of 
(geometric) modules. An $R$-module $M$ over the space ${X}$ is a family $(M_x)_{x \in {X}}$ of finitely 
generated free $R$-modules $M_x$
indexed over points of ${X}$, which is locally finite in the sense 
that $\oplus_{x \in K} M_x$ is finitely generated for every compact subset $K \subset {X}$. A morphism $\phi \colon M \to N$
is an $R$-linear map $\phi \colon \oplus_{x \in {X}} M_x \to \oplus_{y \in X} N_y$. Such a  map can of course be decomposed 
and written as a matrix
$\phi=(\phi_{y,x})$ indexed over ${X} \times {X}$. The additive category of all such modules and morphisms would be denoted 
$\calc({X})$ or $\calc( {X} ; R )$ 
and is equivalent to the category of finitely generated free $R$-modules if ${X}$ is a compact space.

\subsubsection{Support conditions}
We are however only interested in subcategories of modules and morphisms satisfying certain support conditions. 
The support of a module $M$ or a morphism $\phi$ are defined as 
$\supp M = \{ x \in {X} \; | \; M_x \neq 0 \} \subset {X}$ respectively 
$\supp \phi = \{ (x,y) \; | \; \phi_{y,x} \neq 0 \} \subset {X} \times {X}$. For a morphism control condition  $\cale$ 
(a set of subsets of ${X} \times {X}$) and an object support condition
$\calf$ (a set of subsets of ${X}$) we denote by
\[
\calc( {X} , \cale , \calf )
\]
the subcategory of $\calc({X})$ consisting of modules $M$, for which there exists an $F \in \calf$ such 
that $\supp M \subset F$, and morphisms $\phi$ between such modules,
for which there exists an $E \in \cale$ with $\supp \phi \subset E$. 
We will often refer to such morphisms as $\cale$-controlled morphisms. The conditions one needs to impose on $\cale$ and $\calf$ in order to
assure that this yields in fact an additive category are spelled out in Subsection~2.3 of \cite{Bartels-Farrell-Jones-Reich(topology)}. 
A basic example of a morphism control condition on a metric space $({X},d)$ is $\cale_d$, consisting of all $E \subset {X} \times {X}$
for which there is $\alpha > 0$ such that $(x,y) \in E$ implies $d(x,y) < \alpha$.
Measuring control via a map $p \colon {X} \to Y$ is
formalized by pulling back $\cale$ (living on $Y$), i.e.\ forming $p^{-1} \cale = \{ (p \x p)^{-1}( E ) \; | \; E \in \cale \}$.
In the case, where $p$ is the inclusion of a subspace we usually omit $p$ and write
$\cale$ instead of $p^{-1}\cale$. Similar notational conventions apply to the $\calf$'s.

\subsubsection{Equivariant versions}
We usually deal with equivariant versions where ${X}$ is assumed to be a free $\Gamma$-space 
and modules and morphisms are required to be invariant under the $\Gamma$-action. The corresponding category is denoted
\[
\calc^{\Gamma} ( X , \cale , \calf ).
\]
Under suitable conditions about the $\cale$'s and $\calf$'s a $\Gamma$-equivariant map $f \colon {X} \to {Y}$ 
induces a functor on such categories which sends $M$ (a module over ${X}$) to $f_{\ast}M$ (a module over ${Y}$) given by
$f_{\ast} M_y = \oplus_{x \in f^{-1}( \{ y \} )} M_x$.

\subsubsection{Thickenings} \label{subsubsec:thickenings}
If $E$ is a neighborhood of the diagonal in $X \times X$ and $A$ a subset of $X$, then we define the 
$E$-thickening $A^E$ of $A$ in $X$ to be the set of all points $x \in X$
for which there exists a point $a \in A$ such that $(a,x)$ belongs to $E$. 
In the case where $E$ is determined by a constant $\delta$ via a metric or by a pair $(\alpha, \delta)$ 
using the ``foliated distance'' (compare Subsection~\ref{section-fol-control-w-decay}) 
we write $A^{\delta}$, respectively $A^{\alpha, \delta}$ for the corresponding thickenings.

\subsubsection{Germs} \label{subsubsec:germs}
We will often use Karoubi quotients of the categories $\calc^{\Gamma}( X , \cale , \calf )$ 
introduced above: let $\calf_0$ be another object support condition 
which is $\cale$-thickening closed, i.e.\ for every $F \in \calf_0$ and $E \in \cale$ there exists an $F^{\prime} \in \calf_0$
such that $F^E \subset F^{\prime}$. 
Then $\calc^{\Gamma} ( {X} , \cale , \calf )^{> \calf_0}$ is defined as the additive category which 
has the same objects as $\calc^{\Gamma} ( {X} , \cale , \calf )$, but where 
morphisms are identified whenever their difference factors over a module with support in $F_0 \in \calf_0$. 
We think of this construction as taking germs away from $\calf_0$.
If our space is ${X} \x [1,\infty)$ 
and $\calf_0 = \{ {X} \x [1,t] \; | \; t \in [1, \infty) \}$ then we write $\calc^\infty$ rather than $\calc^{> \calf_0}$, 
see Example~\ref{germs-at-infty}. 
Further formal properties of these constructions which will be needed throughout 
the proof are discussed in Appendix~\ref{appendix-fibration-sequences}.

\subsection{Assembly as a ``forget control-map''} \label{subsec:assembly-as-forget-control}

\subsubsection{Resolutions} \label{subsubsec:resolutions}
A $\Gamma$-space is called $\Gamma$-compact if it is the $\Gamma$-orbit
of some compact subspace.
A resolution of the $\Gamma$-space $X$ is a $\Gamma$-map $p \colon \overline{X} \to X$ of $\Gamma$-CW-complexes, 
where $\overline{X}$ is a free $\Gamma$-space and every $\Gamma$-compact set in $X$ is the image under $p$ 
of some $\Gamma$-compact set in $\overline{X}$.  For every space the projection $X \times \Gamma \to X$ is
a functorial resolution called the standard resolution.

\subsubsection{The functor $\cald^{\Gamma}$}
In \cite[Definition~2.7]{Bartels-Farrell-Jones-Reich(topology)} we defined for a not necessarily free $\Gamma$-$CW$ complex $X$ the 
notion of $\Gamma$-equivariant continuous control. This is  a  morphism support condition denoted 
$\cale_{\Gamma cc}(X)$ on the space $X \times [1, \infty)$. 
The set of all $\Gamma$-compact subsets of $\overline{X}$ is denoted $\calf_{\Gamma c}(\overline{X})$. 
The object support condition $p_{\overline{X}}^{-1} \calf_{\Gamma c}(\overline{X})$, where 
$p_{\overline{X}} \colon \overline{X} \times [1, \infty) \to \overline{X}$ denotes 
the projection, is our standard object support condition on $\overline{X} \times [1, \infty)$.
It is shown in 
\cite[Section~3 and 5]{Bartels-Farrell-Jones-Reich(topology)} that up to equivalence 
\[
\cald^{\Gamma}(\overline{X} ; p) = \calc^{\Gamma} ( \overline{X} \times [1, \infty) ; (p \times \id)^{-1} \cale_{\Gamma cc} (X) ,
p_{\overline{X}}^{-1} \calf_{\Gamma c} ( \overline{X} ))
\]
does not depend upon the chosen resolution (hence one can always use $X \times \Gamma \to X$ as a resolution; in this case we denote the
category by $\cald^\Gamma(X)$) and that
$X \mapsto \IK^{-\infty} \cald^{\Gamma} ( \overline{X} ; p)$ yields an equivariant homology theory on the level of homotopy groups.

\subsubsection{Assembly}
The map induced by 
\[
\cald^\Gamma(E \Gamma ( \VCyc )) \to \cald^\Gamma(\pt)
\] 
is on the level of $K$-theory up to an index shift a model for the generalized assembly map that was discussed above.
Compare \cite[Corollary~6.3]{Bartels-Farrell-Jones-Reich(topology)}.


\typeout{---------------------proof -------------------------}

\section{Outline of the proof} \label{section-proof}

The injectivity part of Theorem~\ref{main-theorem} is proven in \cite{Bartels-Farrell-Jones-Reich(topology)}. We will prove surjectivity.
Our first observation is that it suffices to prove surjectivity of the map which is induced on the level of $K$-theory by 
$\cald^\Gamma(X(\infty)) \to \cald^\Gamma(\pt)$, where $X(\infty)$ is any 
$\Gamma$-CW complex all whose nontrivial
isotropy groups are infinite cyclic. In fact by the universal property of $E \Gamma ( \Cyc )$ such a map always  
factorizes over $\cald^\Gamma(E \Gamma ( \Cyc ))$. It follows from \cite[Proposition~3.5]{Bartels-Farrell-Jones-Reich(topology)} 
that instead of considering the map of standard resolutions (compare~\ref{subsubsec:resolutions}) 
we can equally well work with any map of resolutions which covers $X( \infty) \to \pt$.
We proceed to construct a space $X(\infty)$ and such a map of resolutions.

Let $\tilde{M}$ be the universal covering of a closed Riemannian manifold with strictly negative  sectional curvature. 
The hyperbolic enlargement of $\tilde{M}$ is the warped product (compare \cite{Bishop-ONeill(negative-curvature)})
\[
\IH \tilde{M} = \IR \times_{\cosh (t)} \tilde{M}.
\]
It is the differentiable manifold $\IR \times \tilde{M}$ equipped with the Riemannian metric determined by
$ dg^2_{\IH \tilde{M}} = dt^2 + \cosh(t)^2 dg_{\tilde{M}}^2$.
We refer to the $\IR$-factor as the hyperbolic enlargement direction or briefly the $\IH$-direction.
Let $S \IH \tilde M$ 
denote the unit-sphere subbundle of the tangent bundle of $\IH \tilde{M}$. For a subset $A \subset \IR$ we denote
by $S\IH_{A} \tilde{M} $ the restriction of this bundle to the subspace 
$\IH_A \tilde{M} \subset \IH \tilde{M}$ which is defined as $A \times \tilde{M} \subset \IR \times \tilde{M}=\IH \tilde{M}$.
Throughout the paper we also fix the notation
\[
\IB = [ 0 , \infty ) \mbox{ and }  \IT = [ 1 , \infty ).
\]
The space $S \IH \tilde{M} \times \IB \times \IT$ will be important in our context and we will generically use
$h$, $\beta$ and $t$ to denote its $\IH$-, $\IB$- and $\IT$-coordinate.

In Section~\ref{section-collapsing} we factorize the natural projection $S\IH \tilde{M} \to \tilde{M}$
via a map called $p_X$
over a certain $\Gamma$-compact free $\Gamma$-space $X$, i.e.\ we have a commutative diagram
\[
\xymatrix{
S \IH \tilde{M} \ar[r]_-{p_X} \ar@/^3ex/[rr] &
X \ar[r]  & \tilde{M}. \\
         }
\]
Roughly speaking $X$ is obtained from $S \IH \tilde{M}$ by collapsing $S \IH_{(-\infty ,-1]} \tilde{M}$
along the $\IH$-direction to $S \IH_{ \{ -1 \} } \tilde{M}$ and similarly
$S \IH_{[1 , \infty ) } \tilde{M}$ to $S \IH_{\{ 1 \} } \tilde{M}$. Note in particular that 
$S \tilde{M} \subset S \IH_{ \{ 0 \} } \tilde{M}$
sits  naturally as a subspace in $X$.
For details about the map $p_X$ see Subsection~\ref{subsection-construct-X} (resp.\ \cite[14.5]{Bartels-Farrell-Jones-Reich(topology)}).
In Subsection~\ref{subsection-construct-X} we also construct a quotient map 
\[
p\colon X \times \IB \to X( \infty ).
\] 
Here $X( \infty )$ is a $\Gamma$-space all whose isotropy groups are cyclic. 
It is obtained as the infinite mapping telescope of a sequence of maps which collapse more and more lines in $S \tilde{M} \subset X$.
Here the $\IB$-direction is the telescope direction and the lines correspond to preimages of closed 
flow lines of the geodesic flow under the covering projection $S\tilde{M} \to SM$.
(Although $X( \infty)$ is not a model for $E \Gamma ( \Cyc )$ it is fairly close to being one. With more effort one could 
probably work with an actual model at its place.) 
The map of resolutions alluded to above is now
\[
\xymatrix{
X \times \IB \ar[r] \ar[d]^-p & \tilde{ M } \ar[d]^-{\ast} \\
X( \infty ) \ar[r] & \pt .
         }
\]
This map of resolutions
induces the bottom map in the following diagram of additive categories. 
\[ \label{main-diagram}
\xymatrix{
\calc^{\Gamma}( (S \IH \tilde{M} \times \IB \times \IT)_\angle , \cale_w , \calf_{\IB} )^{\infty}
   \ar[ddr]^-{(6)} &
\calc^{\Gamma} ( (S \IH \tilde{M} \times \IT)_\angle , \cale_{geo} )^{\infty} \ar@/^14ex/[dd]^-{(4)}
\ar[l]_-{(5)} \\
\calc^{\Gamma}( (S \IH \tilde{M} \times \IB \times \IT)_\angle , \cale_s , \calf_{\IB} )^{\infty}
\ar[u]_-{(7)}  \ar[d]^-{(8)}  &
\calc^{\Gamma} ( S \IH_{ \{ 0 \} } \tilde{M} 
 \times \IT , \cale_{asy} )^{\infty}  
\ar[u]_-{(3)}  \ar[d]^-{(2)}  \\
\calc^{\Gamma} ( X \times \IB \times \IT , (p \times \id_{\IT})^{-1}
          \cale_{\Gamma cc} (X ( \infty )) , 
            \calf_{\IB} )^{\infty}  \ar[d]^{=} &
\calc^{\Gamma} ( \tilde{M} \times \IT , \cale_d )^{\infty} \ar[d]^-{(1)}
\\
\cald^{\Gamma}( X \times \IB , p ) 
\ar[r]  &
\cald^{\Gamma} ( \tilde{M} , \ast ) .
         }
\]
Here $(S \IH \tilde{M} \x \IT)_\angle$ is the subspace of   $S \IH \tilde{M} \x \IT$
consisting of all points $(v,t)$ where the absolute value of the
$\IH$-coordinate satisfies $|h(v)| \leq t$. Similarly the subspace 
$(S \IH \tilde{M} \x \IB \x \IT)_\angle$ consists of all $(v,\beta,t)$
with $|h(v)| \leq t + \mu_n \beta$.  The constant here is $\mu_n =
10^{n+3}$, where $n$ is the dimension of $S \IH \tilde{M}$.
(The subscript $\angle$ is supposed to remind the reader  of the shape
of the region it describes.) We will mostly be interested in these subspaces, 
but all maps and all object and morphism support conditions 
can and will be defined on the whole spaces. 
The corresponding restrictions to the $\angle$-subspaces will not appear in the notation. 

All maps  in the diagram
except (3) and (8) are induced by the obvious projections, inclusions  or identity maps of the underlying spaces. 
The map (3) is essentially given by the geodesic flow and discussed in detail in Section~\ref{section-flow}.
The map of spaces underlying (8) is induced by $p_X \colon S \IH \tilde{M} \to X$.

The essential information in the diagram is however contained in the different
object and morphism support conditions which will be explained in detail below. 
Going clockwise around the diagram from the lower left hand corner to the lower right hand corner 
should be thought of as forgetting more and more control. 
Our task is to step 
by step gain control going counterclockwise. (The existence of the wrong way maps (5) and (8) says that in between we gain
more control than we actually need.) 

We have the following statements about the diagram above:
\begin{enumerate}
\item \label{dia-commutes}
It will be clear from the construction that the diagram 
without the map (3) commutes.
\item \label{dia-342-h-commutes}
The triangle consisting of the bended arrow (4) and the maps (3) and
(2) induces a triangle in $K$-theory which commutes up to homotopy. 
This will be proven in Corollary~\ref{corollary-triangle-commutes}.
\item \label{dia-2-split-surjective}
The  map (2) induces a split surjection in $K$-theory and hence by (ii) also the maps (4) and (6).
This is an immediate consequence of Proposition~\ref{split-surjective} (compare the discussion before that Proposition).
\item \label{dia-7-equivalence}
The map (7) induces an  equivalence in $K$-theory  by the Foliated Control Theorem~\ref{fct}.
\item \label{dia-easy-1}
The map (1) induces an equivalence 
in $K$-theory by the easy Lemma~\ref{map1}.
\end{enumerate}
These statements  imply that the bottom map induces a
split surjection in $K$-theory and hence our main Theorem~\ref{main-theorem} follows. 

We will now describe the diagram in more detail and explain some
aspects of statements \ref{dia-342-h-commutes}, \ref{dia-2-split-surjective} and 
\ref{dia-7-equivalence}.
We proceed counterclockwise. 
\vspace{1em}

(1) 
Equip $\IT$ with the standard metric. Let $d$ denote any product metric (e.g.\ the $\max$-metric)
on $\tilde{M} \times \IT$ and let $\cale_d$ be the corresponding morphism support condition.
Observe that the continuous-control condition $\ast^{-1} \cale_{\Gamma cc}$ is a weaker condition.
The resulting forget-control map (1) is shown to induce an equivalence in Lemma~\ref{map1} below.
\vspace{1em}

(2)
The base space of the bundle $S\IH_{ \{ 0 \} } \tilde{M}$ is $\tilde{M}$.
The map (2) is induced by the bundle projection.
The space $S \IH \tilde{M}$ comes equipped with two different foliations, the asymptotic foliation $F_{asy}$
and the geodesic foliation $F_{geo}$. These are explained in Section~\ref{section-flow}. In that section we will also
define 
the notion of foliated control with a prescribed decay speed depending on a foliation $F$ and a set of decay speed functions $\cals$. 
The important point about the map (2) is that in the source 
we have foliated control with a certain carefully chosen decay speed (denoted $\cals_{asy}$)
with respect to the asymptotic foliation $F_{asy}$. This kind of control will be denoted $\cale_{asy}$.
Section~\ref{section-transfer} is devoted to proving that the map (2) induces a split surjection in $K$-theory, see in particular
Proposition~\ref{split-surjective}.
To prove this we will construct a transfer map going essentially the other way.
In fact the target of the transfer map will not be $\calc^{\Gamma} ( S \IH_{\{ 0 \} } \tilde{M} \times \IT , \cale_{\asy} )^{\infty}$
but a formally enlarged version of that category which yields the 
same $K$-theory. 
The construction of the transfer map will depend on the choice of a sequence $\Jd =(\delta^0, \delta^1,\delta^2  \dots )$ of 
decay speed functions
which will be responsible for the decay speed one can achieve in the target. The geometry of our situation enters in
verifying that a suitable sequence of decay speed functions exists, compare Lemma~\ref{sequence-exists}. 
The map on the $K$-theory of $\calc^{\Gamma} ( \tilde{M} \times \IT , \cale_d )^{\infty}$
induced by the composition of the transfer and the projection is described by an element in a Swan group in Proposition~\ref{mult-by-swan}.
In order to achieve that this element is  the identity in our case, we have to work with the subbundle $S^+ \IH_{ \{ 0 \} } \tilde{M}$ of  $S\IH_{ \{ 0 \} } \tilde{M}$
whose fiber is contractible. This explains the  necessity of the hyperbolic enlargement: it allows us to pick out this subbundle.

\vspace{1em}

(3)
In Section~\ref{section-flow} we will study the map (3) which is induced by
\begin{eqnarray*}
S \IH \tilde{M} \times \IT & \to & S \IH \tilde{M} \times \IT \\
( v , t ) & \mapsto & ( \Phi_t ( v ) , t ),
\end{eqnarray*}
where $\Phi$ denotes the geodesic flow on $S \IH \tilde{M}$. Via this map one gains control in the directions
transverse to the geodesic flow.
We prove in Theorem~\ref{asytogeo} that the map turns foliated control with respect to the asymptotic foliation 
into foliated control with respect to the geodesic foliation. 
In fact in the target we will have foliated control (with respect to the geodesic foliation) with 
exponential decay speed (depending on the upper bound for 
the sectional curvature) and our choice of decay speed $\cals_{asy}$ for $\cale_{asy}$ in (2) is made in such a way that we achieve this.
Since $\Phi_t ( S \IH_{ \{ 0 \} } \tilde{M} ) \subseteq S \IH_{ [ -t , t ] } \tilde{M} $ 
all objects in the image will lie in $(S \IH \tilde{M} \x \IT)_\angle$.
\vspace{1em}

(4)
The bended arrow (4) is induced by the projection $S \IH \tilde{M} \to \tilde{M}$.
This collapses the non-compact $\IH$-direction and is therefore not a proper map. In general, non proper maps 
do not induce functors on our
categories of modules over a space. 
(Such a map does not preserve the local finiteness condition, compare \ref{subsubsec:modules-and-morphisms}.)  
However, the restriction of the projection to  
$(S \IH \tilde{M} \x \IT)_\angle$ is proper and we obtain a well defined functor on objects.
Since $S \IH \tilde{M} \to \tilde{M}$ does not increase distances and $\cale_{geo}$-controlled morphisms are in particular
bounded with respect to the product metric we obtain a well defined functor.
In Corollary~\ref{corollary-triangle-commutes}
we show that the triangle induced in $K$-theory by the maps (2), (3) and (4) 
commutes up to homotopy. Here we use a variant of the Lipschitz homotopy argument from \cite{Higson-Pedersen-Roe(controlled)}.
\vspace{1em}

(5)
The map (5) is induced by the inclusion of $S \IH \tilde{M} \times \{ 0 \} \times \IT$ into
$S \IH \tilde{M} \times \IB \times \IT$. This inclusion is clearly compatible with the $\angle$-subspaces.
We equip $\IB$ with the standard (Euclidean) Riemannian metric and $S\IH \tilde{M} \times \IB$ with the product Riemannian structure.
Also we extend the geodesic foliation (by taking the product with the trivial $0$-dimensional foliation of $\IB$) to a foliation $F_w$ of $S\IH \tilde{M} \times \IB$.
Now $\cale_w$
is defined similar to $\cale_{geo}$ as foliated control with a certain carefully chosen decay speed $\cals$ (defined in 
Subsection~\ref{subsection-construction-of-S}) with respect to the foliation $F_w$. 
It will follow from the construction (see Proposition~\ref{S-properties}~\ref{S-properties-eins}) 
that the inclusion maps $\cale_{geo}$-control to $\cale_w$-control.  
The $\calf_{\IB}$ object support condition consists of all subsets whose
projection to $\IB$ is contained in a compact interval $[0 , \beta_0]$. 
\vspace{1em}

(6)
The map (6) is induced by the projection $S \IH \tilde{M} \times \IB \to \tilde{ M}$. 
As in (4) this projection is not proper since it collapses the non-compact $\IH$- and
$\IB$-direction. However the restriction to the $\angle$-subspace and
the $\calf_{\IB}$ condition ensure that we nevertheless have a
well defined functor on objects. By construction $\cale_w$-control dominates metric control and the projection does
not increase distances. Therefore (6) is also compatible with the morphism control conditions.
\vspace{1em}

(7)
The fact that the map (7) induces an equivalence in $K$-theory should be thought of 
as a ``foliated control theorem''. It is proven as Theorem~\ref{fct}.
Very roughly, this theorem improves foliated control to metric control (with decay speed) on compact subsets (in the $S \IH \tilde{M}$-coordinate)
that do not meet preimages of ``short'' closed geodesics in $M$.
The relatively long Section~\ref{section-fct} is devoted to the proof of this theorem.
The only difference between source and target of the map (7) are the morphism support conditions.
As explained above the weak control condition $\cale_w$ is essentially, i.e.\ up to the added $\IB$-direction, foliated 
control with respect to 
the geodesic foliation with a certain carefully chosen decay speed $\cals$.
The stronger control condition $\cale_s$ is obtained from $\cale_w$,
by adding a metric control condition with decay speed $\cals$
over  a certain subset $S \subset S \IH \tilde{M} \times \IB \times \IT$. 
The subset and the precise control condition
will be explained in Subsection~\ref{statement-fct}. 
\vspace{1em}

(8)
The map (8) is induced from the projection $p_X \colon S \IH \tilde{M} \to X$. 
In particular this map collapses the non compact $\IH$-direction 
and hence the remaining $\calf_{\IB}$-object support in the target is just the usual
$\Gamma$-compact support. 
The $\cale_s$ condition
is shown to be strong enough to induce a map in Proposition~\ref{pro-strong-maps-to-continous}.
\vspace{1em}

This finishes the outline of the proof.
\vspace{1em}

We start the proof by the following easy lemma about the map (1).

\begin{lemma} \label{map1}
The map (1) induces an equivalence in $K$-theory. 
\end{lemma}
\begin{proof} The map $\calc^{\Gamma}( \tilde{M} \times \IT , \cale_d ) \to 
\calc^{\Gamma} ( \tilde{M} \times \IT , \ast^{-1} \cale_{\Gamma cc} )$ induces a map between the two
corresponding ``germs~at~infinity''-fibrations, see Example~\ref{germs-at-infty}. In the resulting ladder the two middle terms allow
an Eilenberg-swindle towards infinity along $\IT$ (compare \cite[4.4, 4.5]{Bartels-Farrell-Jones-Reich(topology)}) and the left hand terms are even equal.
\end{proof}


\typeout{---------------------flow -------------------------}

\section{Gaining control via the geodesic flow} 
\label{section-flow}

There are two foliations on $S \IH \tilde{M}$, the geodesic foliation and the asymptotic foliation. 
In this section we will define the notion of foliated control with a prescribed decay speed and we will
prove in Theorem~\ref{asytogeo} that the geodesic flow can be used to turn foliated control with a certain decay speed 
with respect to the asymptotic foliation into foliated control with exponential decay speed with respect to the 
geodesic foliation. Finally we show that after forgetting  control the map induced by the geodesic flow is homotopic 
to the identity (more precisely to a certain inclusion), see Theorem~\ref{lipschitz}.

\subsection{Geometric preparations}

Recall from Section~\ref{section-proof} that $\IH \tilde{M}$ denotes the hyperbolic enlargement 
and $S \IH \tilde{M}$ its sphere bundle.
The space $S \IH \tilde{M}$ will be equipped with two foliations, the geodesic foliation and the 
asymptotic foliation. 
Let $\Phi: \IR \times S \IH \tilde{M} \to S \IH \tilde{M}$, $(t,v) \mapsto \Phi_t(v)$ denote the geodesic flow. 
The geodesic foliation
$F_{geo}$ is simply the $1$-dimensional foliation given by the flow lines of $\Phi$.
Two points $v$ and $w$ in $S \IH \tilde{M}$ are 
called asymptotic if the distance between $\Phi_t (v)$ and $\Phi_t (w)$ stays bounded if $t$ tends to $+ \infty$.
This defines an equivalence relation and 
the set of equivalence classes, denoted $S(\infty)$, can be naturally equipped with a topology in such a way
that the map $a: S \IH \tilde{M} \to S( \infty)$ given by sending a vector to its equivalence class restricts to
a homeomorphism on each fiber $S \IH \tilde{M}_x$ of the bundle $S \IH \tilde{M} \to \IH \tilde{M}$, compare Section~1
in \cite{Eberlein-ONeill(visibility)}. The preimages $a^{-1}(\theta)$ for $\theta \in S( \infty )$
are the leaves of a foliation $F_{asy}$ which we will call the asymptotic foliation. 
Since $M$ is compact there are positive constants $a$ and $b$ such that the sectional curvature $K$ satisfies 
\begin{eqnarray} \label{curvature-constants}
-b^2 \leq K \leq -a^2.
\end{eqnarray}
The same inequalities hold for the hyperbolic enlargement (compare \cite{Farrell-Jones(dynamics-I)}). 

The homeomorphisms $S \IH \tilde{M}_x \to S( \infty )$  
are used to define the fiber transport 
\[
\nabla_{y,x} : S \IH \tilde{M}_x \xrightarrow{\cong} S( \infty ) \xleftarrow{\cong}  S \IH \tilde{M}_y
\]
for the bundle $S \IH \tilde{M} \to \IH \tilde{M}$.
The fiber transport will play an important role in Section~\ref{section-transfer}. 
Since we have curvature bounds the fiber transport
is known to be H{\"o}lder-continuous. More precisely a consequence of  
Proposition~2.1.\ in \cite{Anderson-Schoen(positive-harmonic)} is the following lemma which will be used in Theorem~\ref{asytogeo} and Lemma~\ref{sequence-exists}.

\begin{lemma} \label{hoelder}
For all $\alpha > 0$ there is a constant $C_0(\alpha) > 0$ 
such that for $x$, $y \in \IH \tilde{M}$ with $d(x,y) \leq \alpha$ 
and $v$, $w \in S \IH \tilde{M}_x$
we have
\[
d( \nabla_{y,x} (v), \nabla_{y,x} (w) ) \leq C_0(\alpha) \cdot d(v,w)^{\frac{a}{b}} .
\]
\end{lemma}

Later we will have to quantitatively analyze how the flow $\Phi_t:S\IH \tilde{M} \to S \IH \tilde{M}$ may increase 
distances. 
For this purpose we introduce the 
function $C_{flw} (t)$ in the following lemma.

\begin{lemma} \label{flowconstant}
There exists a monotone increasing function $C_{flw}(t)$ such that $|d\Phi_s| \leq C_{flw}(t)$ for all 
$|s| \leq t$. In particular for arbitrary $v$, $w \in S \IH \tilde{M}$ we have
\[
d(\Phi_t(v), \Phi_t(w)) \leq C_{flw}(t)  \cdot d(v,w) .
\]
\end{lemma}

\begin{proof}
This again holds since we have curvature bounds.
The differential of the geodesic flow can be  expressed 
in terms of Jacobi fields, compare \cite[Section 2.3]{Eberlein-Hamenstaedt-Schroeder(nonpositive)}. 
These satisfy a second order differential equation involving the sectional curvature as 
coefficients \cite[p.15-16]{Cheeger-Ebin(comparison)} and the result can be deduced from this equation using standard 
arguments about ordinary differential equations, compare e.g.\ \cite[p.79]{Perko(differential)}.
\end{proof}

\subsection{Foliated control with decay speed}
\label{section-fol-control-w-decay}

We now want to define the notions of metric respectively foliated control with decay speed $\cals$. 
Let $F$ be a foliation of a Riemannian manifold $N$. For $x,y \in N$ we will write  
\[
d_F(x,y) \leq (\alpha,\delta)
\]
if there is a piecewise $C^1$-path 
of arclength  shorter than $\alpha$ which is entirely contained in one leaf of the 
foliation and whose start- respectively end-point lies  within distance $\frac{\delta}{2}$ of $x$ respectively $y$. 
(Elsewhere we used $\delta$ instead of $\frac{\delta}{2}$ but compare Remark~\ref{alpha-zero-control}.)
Suppose we are given a set 
$\cals$ of functions from $\IT$ to $[0,\infty)$. (We often use $\delta_t$ as the name for the function which sends $t$ to $\delta_t$.) 
Suppose $\cals$ satisfies the following conditions.
\begin{enumerate} \label{AandB}
\item[(A)]
For each $\delta_t \in \cals$ and every $\alpha \in \IR$ there exists $\delta_t^{\prime}$ and $t_0 \geq 1$ such 
that $\delta_{t + \alpha} \leq \delta_t^{\prime}$ for all $t \geq t_0 + |\alpha|$.
\item[(B)]
Given $\delta_t$, $\delta_t^{\prime} \in \cals$ there exists $\delta_t^{\prime \prime} \in \cals$ and $t_0 \geq 1$
such that $\delta_t + \delta_t^{\prime} \leq \delta_t^{\prime \prime}$ for all $t \geq t_0$.
\end{enumerate}
Given such a set of functions $\cals$ we make the following definitions.

\begin{definition}[Foliated and metric control with decay speed] 
\label{fol-control-w-decay}
$\mbox{ }$
\begin{enumerate}
\item
Suppose $X=(X,d)$ is a metric space, then we define a morphism support condition 
$\cale=\cale(X,\cals)$ on $X \times \IT$ by requiring
that a subset $E$ of $(X \times \IT)^{\times 2}$ belongs to $\cale$ if there exists a  
$\delta_t \in \cals$ and  constants $\alpha>0$ and $t_0 > 1$ such that for all $(x,t,x^{\prime},t^{\prime}) \in E$ 
we have 
\[
|t-t^{\prime}| \leq \alpha, \quad d(x,x^{\prime}) \leq \alpha
\]
and if $t,t' > t_0$ then
\[
d(x,x^{\prime} ) \leq \delta_{\min(t,t^{\prime})} .
\]
This condition will be called {\em  metric control with decay speed} $\cals$.
\item
Suppose $N$ is a Riemannian manifold equipped with a foliation $F$. 
We define $\cale=\cale(N,F,\cals)$ a set of subsets of $(N \times \IT)^{\times 2}$ by requiring that a subset 
$E$ belongs to $\cale$ if there exists a function $\delta_t \in \cals$ and 
constants $\alpha>0$ and $t_0 > 1$ such that for all
$(x,t,x^{\prime},t^{\prime}) \in E$ we have 
\[
|t-t^{\prime}|\leq \alpha, \quad  d(x,x') \leq \alpha
\]
and 
if $t,t' > t_0$ then
\[
d_{F} (x,x^{\prime} ) \leq (\alpha , \delta_{\min(t,t^{\prime})}) .
\]
If $\cale$ defines a morphism-control condition then it will be called
{\em foliated control with respect to the foliation} $F$ \em{with decay speed} $\cals$.
\end{enumerate}
\end{definition}

\begin{remark} \label{alpha-zero-control}
Observe that for $\alpha=0$ the foliated condition in (ii) specializes to the metric condition in (i), 
i.e.\ if $\alpha=0$ then an $(\alpha, \delta_t)$-foliated controlled morphism is $\delta_t$-controlled in the 
metric sense.
\end{remark}

\begin{remark} \label{only-need-germs}
Note that these definitions  only depend on the behaviour of functions in $\cals$ in a neighborhood of $\infty$, 
i.e\ what is  really important about $\cals$
is the set of germs (at infinity) of functions it determines. 
\end{remark}

\begin{warning} \label{warning-fol-control}
While conditions (A) and (B) guarantee that in the metric case $\cale(X,\cals)$ is a morphism-control condition, in general 
this may fail in the foliated case for $\cale = \cale(N,F,\cals)$ because it is not 
clear that $\cale$ is closed under composition.
However in the two cases we are interested in $\cale$ is closed under composition by Lemma~\ref{foltriangle} combined 
with Lemma~\ref{construct-S-asy}~\ref{construct-S-asy-2} below. 
\end{warning}

\subsection{Gaining control via the geodesic flow}

Now set
\begin{eqnarray*}
\cals_{geo} & = & \{ A \cdot \exp (-at)  |   A > 0 \}, \\
\cals_{asy} & = & \{ A \cdot \exp ( -\lambda ( C_{flw}(2t+B) + (t+B)^2 )) | A > 0, B \in \IR, \lambda > 0 \} ,
\end{eqnarray*}
where in the first line $a$ comes from the upper curvature bound in (\ref{curvature-constants}).
In the second line we take  $C_{flw}(2t+B) = C_{flw}(0)$ for $2t+B<0$. 
Clearly \label{no comma!} conditions (A) and (B) before
Definition~\ref{fol-control-w-decay} are satisfied in both cases.
Here $\cals_{asy}$ is designed in such a way that the geodesic flow will turn $\cals_{asy}$-decay speed
into $\cals_{geo}$-decay speed, see Theorem~\ref{asytogeo} which uses
property \ref{construct-S-asy-1} of the following elementary lemma.

\begin{lemma}  \label{construct-S-asy}
Apart from (A) and (B) before Definition~\ref{fol-control-w-decay} $\cals_{asy}$ satisfies
\begin{enumerate}
\item \label{construct-S-asy-1}
For $\delta_t \in \cals_{asy}$ and each $\alpha > 0$ there exists a $\delta_t^{\prime} \in \cals_{geo}$ and $t_0 \geq 1$
such that 
\[
C_{flw}(t+ \alpha) \cdot \delta_t \leq \delta_t^{\prime}  \mbox{ for all $t \geq t_0$ }.
\]
\item \label{construct-S-asy-2}
If $\delta_t \in \cals_{asy}$ and $\lambda > 0$ then $(\delta_t)^{\lambda} \in \cals_{asy}$.
\end{enumerate}
\end{lemma}

We now define $\cale_{geo}$ to be foliated control with exponential decay speed $\cals_{geo}$ 
with respect to the geodesic foliation $F_{geo}$
and $\cale_{asy}$ to be foliated control with decay speed $\cals_{asy}$ with respect to the asymptotic foliation 
$F_{asy}$, i.e. 
\begin{eqnarray*}
\cale_{geo} & = & \cale (S \IH \tilde{M} , F_{geo} , \cals_{geo} ),    \\
\cale_{asy} & = & \cale (S \IH \tilde{M} , F_{asy} , \cals_{asy} ).   
\end{eqnarray*}
We are now prepared to formulate the main result of this section. 
Recall that $(S \IH \tilde{M} \x \IT)_\angle$ is the subspace given by
$|h| \leq t$, 
where $h$ denotes the $\IH$-coordinate
and $t$ the $\IT$-coordinate. 

\begin{theorem} \label{asytogeo}
The map $(v , t) \mapsto (\Phi_t (v) , t )$ on $S \IH \tilde{M} \times \IT$ turns $\cale_{asy}$-control into
$\cale_{geo}$-control. In particular, it induces a well defined map
\[
\calc^{\Gamma} ( S \IH_{ \{ 0 \} } \tilde{M} \times \IT , \cale_{asy} )^{\infty} \to
\calc^{\Gamma} ( (S \IH \tilde{M} \times \IT)_\angle , \cale_{geo})^{\infty} .
\]
\end{theorem}

\begin{proof}
We recall some results which were discussed in 
\cite{Bartels-Farrell-Jones-Reich(topology)} in
Proposition~14.2 and Lemma~14.3 and rely on \cite{Heintze-ImHof(horospheres)}.
Let $a$ and $b$ be the curvature constants from (\ref{curvature-constants}). With the constants $C=(1+b^2)^{1/2}$, 
$D=1/a$ and the function $E(\alpha)=\frac{2}{b} \sinh ( \frac{b}{2} (\alpha + \frac{1}{a} ) )$ we have for any 
pair $v$, $w \in S \IH \tilde{M}$ of  asymptotic vectors with $d(v,w)\leq \alpha$ the following inequality
\[
d_{F_{geo}} ( \Phi_t ( v ) , \Phi_{t^{\prime}} (w) ) \leq 
\left( C \cdot ( \alpha + |t-t^{\prime}| + D ), C \cdot 2 E( \alpha ) \cdot e^{-at} \right).
\]
Now consider $(v,t,w,t^{\prime})$, where $v$ and $w$ are no longer assumed to be asymptotic.
Assume $t^{\prime} > t$ then we have because of Lemma~\ref{flowconstant} for any monotone decreasing function 
$\delta_t$ that 
\[
d_{F_{asy}} (v , w ) \leq ( \alpha , \delta_t ) \mbox{ and } t^{\prime} -t \leq \alpha 
\]
implies
\[
d_{F_{geo}} ( \Phi_t(v), \Phi_{t^{\prime}} ( w ) ) \leq 
\left( C \cdot ( 2 \alpha + D ) , C \cdot 2 E(\alpha ) \cdot e^{-at} + C_{flw}(t + \alpha ) \delta_t \right).
\]
By Lemma~\ref{construct-S-asy}~\ref{construct-S-asy-1} and property (B) for $\cals_{geo}$ this implies the claim about the morphism control conditions.
Since the flow-speed in the $\IH$-direction is at most $1$ we see that
$(v , t) \mapsto (\Phi_t (v) , t )$ maps $S \IH_{ \{ 0 \} } \tilde{M} \times \IT$ to 
$(S\IH\tilde{M} \x \IT)_\angle$.
\end{proof}

We still have to verify that $\cale_{geo}$ and $\cale_{asy}$ are well defined morphism-support conditions, compare
Warning~\ref{warning-fol-control}.
This is a consequence of Lemma~\ref{construct-S-asy} and  the following lemma.

\begin{lemma}[Foliated triangle inequalities] \label{foltriangle}
For $u$, $v$, $w \in S \IH \tilde{M}$ we have:
\begin{enumerate}
\item \label{foltriangle-eins}
If $d_{F_{geo}}(u,v) \leq (\alpha,\delta)$ and $d_{F_{geo}}(v,w) \leq (\beta,\epsilon)$ then
\[
d_{F_{geo}}(u,w) \leq (\alpha + \beta, \delta +  C_{flw}(\alpha)({\epsilon + \delta}) + \epsilon).
\]
\item
If $d_{F_{asy}}(u,v) \leq (\alpha,\delta)$ and $d_{F_{asy}}(v,w) \leq (\beta,\epsilon)$ then
\[
d_{F_{asy}}(u,w) \leq (C \alpha + C(\frac{\delta + \e}{2}) + \beta, 
         \delta + 2 C_0(\alpha + \frac{\delta + \e}{2}) (C+1)^{\frac{a}{b}} 
                                   (\frac{\delta+\e}{2})^{\frac{a}{b}} + \epsilon).
\]
\end{enumerate}
Here the constant $C_0$ stems from Lemma~\ref{hoelder} and $C$ from the proof of Theorem~\ref{asytogeo}.
\end{lemma}

\begin{proof}
In both cases there are $u_0$, $v_0$, $v_1$ and $w_1$ within distance $\delta$ respectively $\e$ from $u$, $v$ 
respectively $w$ such that $u_0$ and $v_0$ respectively $v_1$ and $w_1$ can be joined by a curve of length $\alpha$ 
respectively $\beta$ contained in a leaf of the foliation in  question. For the first statement assume that 
$\Phi_t(u_0) = v_0$ and let $\bar{u}_0 = \Phi_{-t}(v_1)$. Then $\bar{u}_0$ and $w_1$ are
contained in a leaf of $F_{geo}$ and can be used to  prove the first inequality (using \ref{flowconstant}). 
For the second statement let $x$ and $y$ 
denote the  foot points of $u_0$ and $v_1$ and set $\bar{u}_0 = \nabla_{x,y}(v_1)$. 
Then $\bar{u}_0$ and $w_1$ are
contained in a leaf of $F_{asy}$ and can be used to  prove the second inequality (using 
\ref{hoelder} and \cite[14.3]{Bartels-Farrell-Jones-Reich(topology)}). 
\end{proof}

We finish this section by comparing the map induced by the flow to the map induced by the inclusion
$S \IH_{ \{ 0 \} } \tilde{M} \to S \IH \tilde{M}$.
This is only possible after relaxing the control conditions in the target.
We denote by $\cale_d$ the control condition coming from the product metric on $S \IH \tilde{M} \times \IT$.

\begin{theorem}   \label{lipschitz}
The map $S \IH_{ \{ 0 \} } \tilde{M} \times \IT \to (S \IH \tilde{M} \times \IT)_\angle$ defined
by $(v,t) \mapsto (\Phi_t(v),t)$ and the inclusion induce homotopic maps 
\[
\calc^{\Gamma} ( S \IH_{ \{ 0 \} } \tilde{M} \times \IT , \cale_{asy} )^{\infty} \to
\calc^{\Gamma} ( (S \IH \tilde{M} \times \IT)_\angle , \cale_d )^{\infty} .
\]   
on the level of $K$-theory.
\end{theorem}

\begin{proof}
Let us abbreviate the two categories from above by $\calc_0$ and $\calc_1$. 
Let $\hat Z \subset S \IH _{ \{ 0 \} } \tilde{M} \times \einsu \times \IT$ 
consist of all $(v,s,t)$ with $s \leq t$. Let $p : \hat{Z} \to S \IH _{ \{ 0 \} } \tilde{M} \times \IT$
denote the obvious projection.
We will use 
$\hat\calc = \calc^{\Gamma} (  \hat Z,  \hat\cale )^{\infty}$.
Here $\hat\cale = \cale_d \cap p^{/1} \cale_{asy}$ where $\cale_d$ denotes metric control with respect to
a product metric on $\hat{Z}$. 
The arguments used in the proof of Theorem~\ref{asytogeo} can also be 
used to check that $(v,s,t) \mapsto (\Phi_s(v),t)$ induces a functor $H : \hat\calc \to \calc_1$. 
Moreover, $(v,t) \mapsto (v,1,t)$ and $(v,t) \mapsto (v,t,t)$ induce functors $I,J : \calc_0 \to \hat\calc$ while
$(v,s,t) \mapsto (v,t)$ induces $P : \hat\calc \to \calc_0$. The claim of the theorem is that $H \circ I$ and 
$H \circ J$ induce equivalent maps in $K$-theory. Clearly, $P \circ I = P \circ J = \id_{\calc_0}$. 
It is now sufficient to show that $I$ induces an isomorphism in $K$-theory, since then $I \circ P = \id_{\hat\calc}$ and
$H \circ I = H \circ I \circ P \circ J = H \circ J$ in $K$-theory.
Now $I : \calc_0 \to \hat\calc$ is equivalent to a Karoubi filtration with quotient
$\calc^\Gamma( \hat Z, \hat\cale )^{> \calf_q}$, 
where $\calf_q$ consists of all sets of the form 
$\{ (v,s,t) | s \leq t \mbox{ and } ( s \leq N \mbox{ or } t \leq N ) \}$
for some $N$. We claim that this category is flasque. Indeed, the map $(v,s,t) \mapsto (v,s-1,t)$ 
induces an Eilenberg swindle on it.  
A little care is needed in producing the swindle from this map, since it not well-defined on  $(v,s,t)$ for $s < 2$ 
(because then $s-1 \not\in \einsu$). However, in the quotient category in question modules over 
$S \IH _{ \{ 0 \} } \tilde{M} \times [1,2] \times \IT$ can be ignored. 
Compare~\ref{subsubsec:germs}.
\end{proof}

\begin{corollary} \label{342-commutes}
The triangle consisting of the maps (3),(4) and (2) in the main diagram of Section~\ref{section-proof}
commutes up to homotopy after applying $K$-theory. 
\end{corollary}

\begin{proof}
Compose the maps in Theorem~\ref{lipschitz} with the map induced by the projection 
$S \IH \tilde{M} \x \IT \to \tilde{M} \x \IT$. 
\end{proof}


\typeout{---------------------transfer -------------------------}

\section{The transfer} \label{section-transfer}

Our aim in this section is to prove that the map (2) in our main diagram in Section~\ref{section-proof}
induces a 
split surjective map in K-theory. We define $S^+ \IH \tilde{M}$ to be the subbundle of the sphere 
bundle $S \IH \tilde{M} \subset T \IH \tilde{M} =  T \IR \times T \tilde{M} = \IR \times T \tilde{M}$ consisting of all 
vectors with non-negative $\IR$-coordinate. Note that the fiber of this subbundle is a disk and hence contractible. 
This is important, because we will show below that the transfer on a bundle whose fiber has interesting topology 
is in general not a splitting of the bundle projection, cf.\ Proposition~\ref{mult-by-swan}. 
Since the projection  $S^+ \IH_{ \{ 0 \} } \tilde{M} \to \tilde{M}$ factorizes 
as 
\[
\xymatrix{
S^+ \IH_{ \{ 0 \} } \tilde{M} \ar@{^{(}->}[r] &   S \IH_{ \{ 0 \} } \tilde{M}  \ar[r] &  \tilde{M}
         }
\]
surjectivity of the map (2) is implied by
the following proposition.

\begin{proposition} \label{split-surjective}
The map
\[
\calc^{\Gamma} ( S^+ \IH_{ \{ 0 \} } \tilde{M}  \times \IT , \cale_{asy} )^{\infty} 
\to \calc^{\Gamma} (\tilde{M} \times \IT ; \cale_d )^{\infty} 
\]
induced by the bundle projection
induces a split surjective map in $K$-theory. 
\end{proposition}

In order to prove this proposition we will produce a transfer map in the reverse direction. In fact we will construct 
the following (non-commutative) diagram.
\begin{eqnarray} \label{diagram-transfer}
\xymatrix{
\calc^{\Gamma} ( S^+ \IH_{ \{ 0 \} } \tilde{M}  \times \IT , \cale_{asy} )^{\infty} \ar[r]  \ar[d] & 
\widetilde{ \ch}_{\hf}  \calc^{\Gamma} ( S^+ \IH_{ \{ 0 \} } \tilde{M}  \times \IT , \cale_{asy} )^{\infty} \ar[d] \\
\calc^{\Gamma} (\tilde{M} \times \IT ; \cale_d )^{\infty} \ar[r] \ar[ur]^-{\tr^{\Jd}}  &
\widetilde{ \ch}_{\hf}  \calc^{\Gamma} (\tilde{M} \times \IT ; \cale_d )^{\infty} .
         }
\end{eqnarray}
The diagonal arrow is the promised transfer. It will depend upon the choice of a sequence 
$\Jd=( \delta^0, \delta^1 , \delta^2 , \dots )$ of decay speed function from $\cals_{\asy}$. The horizontal arrows 
induce equivalences in $K$-theory  and the square without the transfer 
map commutes. In Subsection~\ref{subsection-mult-with-swan} we will show that the lower triangle commutes in $K$-theory
up to multiplication by a certain element in a Swan group which is determined by the homology of the fiber of the 
bundle $S^+ \IH_{ \{ 0 \} } \tilde{M}$. Since the fiber of this bundle is a disk we know that the triangle induces a 
commutative triangle, see Corollary~\ref{corollary-triangle-commutes}.

\begin{remark} \label{remark-only-connected}
In the diagram above and in the proof below we have to deal with certain Waldhausen categories (categories with 
cofibrations and weak equivalences \cite{Waldhausen(1126)}) which are categories of 
chain complexes. There seems to be no good definition of non-connective $K$-theory in the literature 
which applies in this generality. 
However in our situation we can always make an ad-hoc construction of a non-connective $K$-theory spectrum as follows. 
Let $X$ be a free $\Gamma$-space and $\cale$ a control condition on $X$. Let $\cale_d$ be the standard Euclidean 
metric control condition on 
$\IR^n$. Let $p: \IR^n \times X \to \IR^n$ and $q: \IR^n \times X \to X$ be the projections and $\IK( - )$ be 
Waldhausen's connective $K$-theory functor which applies to Waldhausen categories. It is well known that the spaces 
\[
\IK \calc^{\Gamma} ( \IR^n \times X ; p^{-1} \cale_d \cap p^{-1} \cale )
\]
together with structure maps derived from swindles coming from a decomposition $\IR^n = \IR^n_+ \cup \IR^n_-$
yield a model for the non connective $K$-theory spectrum of $\calc^{\Gamma} ( X ; \cale )$ (compare the last page in 
\cite{Cardenas-Pedersen(karoubi)}). The same construction applies to categories like $\ch_f \calc ( X ; \cale )$ and 
$\ch_{\hf} \calc ( X ; \cale)$ and all variants which will be used below.  Hence in each case it makes sense to talk 
of the non-connective $K$-theory. In all constructions and arguments below the $\IR^n$-factor will play the role of a 
dummy variable. In order to facilitate the exposition we hence formulated all arguments only for the $0$-th spaces, 
i.e.\ for connective $K$-theory. It is straightforward to make the necessary modifications to obtain the analogous 
statement for the other spaces of the spectrum and to check compatibility with the structure maps.
\end{remark}

\subsection{Set-up}

The construction of the transfer works in the following generality. Let $\tilde{B}$ be the universal covering of the 
compact space $B$. Let $\Gamma$ denote the fundamental group which acts on $\tilde{B}$ from the left. 
Let $\pi: \tilde{B} \to B$ denote the covering projection. Suppose $p:E \to B$ is a smooth fiber bundle with compact 
fiber. We form the pullback $\tilde{E}$ and use the following notation.
\[
\xymatrix{
\tilde{E} \ar[d]^-{\tilde{p}} \ar[r]^-{\pi} & E \ar[d]^-{p} \\
\tilde{B} \ar[r]^-{\pi} & B.
         }
\]
Suppose $\tilde{B}$ is equipped with a $\Gamma$-invariant metric $d$.
Suppose $E$ is a Riemannian manifold and $\tilde{E}$ is equipped with the pulled back Riemannian structure
and a $\Gamma$-invariant foliation $F$. Furthermore let $\cals$ be a set of decay speed functions 
(compare Subsection~\ref{section-fol-control-w-decay}) and assume that $\cale=\cale(\tilde{E},F,\cals)$ really defines a 
morphism control condition, compare Warning~\ref{warning-fol-control}.

We also assume that we are given a fiber transport $\nabla$, i.e.\ a homeomorphism of fibers 
$\nabla_{b',b}: \tilde{E}_b \to \tilde{E}_{b'}$ for each pair of points $b$ and $b'$ in $\tilde{B}$ which fulfills 
the following requirements.

\begin{assumption}  \label{assumption-on-nabla}
The fiber transport has the following properties:
\begin{enumerate}
\item \label{assumption-functoriality}
It is functorial, i.e.\ $\nabla_{b,b}=\id_{\tilde{E}_b}$ and $\nabla_{b'',b'} \circ \nabla_{b',b} = \nabla_{b'',b}$ for all 
$b''$, $b'$ and $b \in \tilde{B}$.
\item
It is $\Gamma$-invariant, i.e.\ for all $b$, $b' \in \tilde{B}$ and all $g \in \Gamma$ we have 
\[
l_g \circ \nabla_{b',b} \circ l_g^{-1} = \nabla_{g b' , gb}.
\]
Here $l_g: \tilde{E}_b \to \tilde{E}_{gb}$ is the restriction of the left action of 
$\Gamma$ on $\tilde{E}$.
\item \label{assumption-transport-metric}
It is compatible with the foliation in the following strong sense. There exists a constant $C\geq 1$ such that 
for all $b,b' \in \tilde{B}$ with $d(b,b') \leq \alpha$ and every $e \in \tilde{E}_b$ we have 
\[
d_{\calf} ( \nabla_{b',b} e , e ) \leq ( C \alpha , 0 ),
\]
i.e.\ there is a path of length no longer than $C \alpha$ inside one leaf which connects $\nabla_{b , b'} e$ and $e$. 
\end{enumerate}
\end{assumption}

An additional requirement will be formulated in Assumption~\ref{assumption-about-d} below.
All the assumptions are fulfilled in our situation where $E \to B$ is the bundle $S^+ \IH_{ \{ 0 \} } M \to M$, the foliation 
is the asymptotic foliation $F_{\asy}$ and $\cals= \cals_{\asy}$. Assumption~\ref{assumption-transport-metric} 
follows from \cite[1.1]{Ballmann-Brin-Eberlein(nonpositive)}: if $\phi(t)$ is a geodesic in $\tilde{M}$ from $b$ to $b'$,
then the path $t \mapsto \nabla_{\phi(t),b} e$ is contained in a leaf of $F_{asy}$ and no longer than 
$(1+b^2)^{1/2} \cdot \alpha$, where $-b^2$ is a lower bound for the curvature of $M$.

\subsection{Homotopy finite chain complexes} \label{subsection-homotopy-finite-chain-complexes}

Below we would like to work with singular chain complexes, which have the advantage that they do not depend on
any further choices (like triangulations or CW-structures). On the other hand in order to define 
$K$-theory we need to impose some finiteness conditions, i.e.\ we want to work with homotopy finite
chain complexes. We now introduce the necessary notation.
Given a free $\Gamma$-space $X$ and a control-condition $\cale$ we define
$\overline{\calc}^{\Gamma} ( X ; \cale )$ completely analogous to $\calc^{\Gamma} ( X ; \cale )$ but we do not
require that the modules are locally finite. We hence allow objects $M=(M_x)$ whose support 
$\supp M = \{ x \in X \; | \; M_x \neq 0 \}$ is an arbitrary subset of $X$. Also the free $R$-modules $M_x$
need not be finitely generated.
(We should however require that the cardinality of the bases are bounded by some fixed large enough cardinal.
This allows us to choose small models for all the categories that will appear below.)
We think of the full subcategory $\calc^{\Gamma} ( X ; \cale )$ inside  $\overline{\calc}^{\Gamma} ( X ; \cale )$ 
as the category of ``finite'' objects and define as explained in the Appendix~\ref{appendix-homotopy-finite}
the categories of finite, respectively homotopy finite chain complexes
\[
\ch_f  \calc^{\Gamma} ( X ; \cale ) \quad \mbox{ and } \quad \ch_{\hf} \calc^{\Gamma} ( X ; \cale ) .
\]
Both categories are full subcategories of the category
$\ch \overline{\calc}^{\Gamma} ( X ; \cale )$
of ``all'' chain complexes and are naturally equipped with the structure of a Waldhausen category, 
see Appendix~\ref{appendix-homotopy-finite}.
The natural inclusions
\[
\calc^{\Gamma} ( X ; \cale )  \to \ch_f \calc^{\Gamma} ( X ; \cale ) \to \ch_{\hf} \calc^{\Gamma} ( X ; \cale ) 
\]
induce equivalences on $K$-theory, compare Lemma~\ref{lemma-inclusions-equivalences} and Remark~\ref{remark-only-connected}.
Analogous considerations apply to $\calc^{\Gamma} ( X \times \IT ; \cale )^{\infty}$ considered as the subcategory 
of ``finite'' objects in $\overline{\calc}^{\Gamma} ( X \times \IT ; \cale )^{\infty}$.

\subsection{The fiber complex}

A chain complex $C \in \ch \overline{\calc}^{\Gamma} ( \tilde{E} \times \IT , \cale )$ is called a fiberwise
chain complex if no differential connects different fibers, i.e.\ if a pair of points  lies in the support
of a differential, then both points lie in the same fiber of the bundle $\tilde{E} \times \IT \to \tilde{B} \times \IT$. 
Given such a fiberwise complex and a point $(b,t)$ in the base we define the fiber $C_{(b,t)}$ to be the largest 
subcomplex such that the support of all its modules lie in the fiber 
$\tilde{E} \times \IT_{(b,t)} = \tilde{E}_b \times \{ t \}$. 

We define the fiberwise complex $F$ and for a given $\delta = \delta_t \in \cals$ the fiberwise complex $F^{\delta}$ by
\[
F_{(b,t)}  =  C_{\sing} ( E_{\pi(b)} ) \quad \mbox{ and } \quad
F^{\delta}_{(b,t)}  =  C_{sing}^{\delta_t}( E_{\pi(b)} ).
\]
Here $C_{\sing}$ denotes the singular chain complex and $C_{\sing}^{\delta_t}$ denotes the subcomplex generated by all 
singular simplices $\sigma: \Delta \to E_{\pi(b)} $ which have the property that  
$\sigma ( \Delta )$ has diameter $\leq \delta_t$ in $E$. \label{AENDERUNG!}
The complexes $F_{(b,t)}$ and $F_{(b,t)}^{\delta}$ are 
complexes  over $E_{\pi(b)}$ by gluing each singular simplex to the image of its barycenter. The complexes $F$ and $F^{\delta}$
become $\Gamma$-invariant complexes over $\tilde{E} \times \IT$ via the maps
\[
\xymatrix{
E_{\pi(b)} & \ar[l]_-{\pi}^-{\cong} \tilde{E}_b \ar[r]^-{\cong} & \tilde{E}_b 
                                  \times \{ t \} \ar[r]^-{\inc} & \tilde{E} \times \IT
         }.
\]
Observe that for $F^{\delta}$ the condition on the size of the singular simplices assures that each differential is 
$(0, \delta_t)$- and hence $\cale=\cale( F , \cals )$-controlled so that 
$F^{\delta} \in \ch \overline{\calc}^{\Gamma} ( \tilde{E} , \cale )$. This is not true for the full singular chain 
complex $F$.

Given an $R$-module $M \in \calc^{\Gamma} ( \tilde{B} \times \IT , \cale_d ; R)$ and a fiberwise $\IZ$-chain complex 
$C \in \ch \overline{\calc}^{\Gamma} ( \tilde{E} \times \IT , \cale ; \IZ)$ we define the fiberwise $R$-chain complex 
\[
M \otimes C \in \ch \overline{\calc}^{\Gamma} ( \tilde{E} \times \IT , \cale ; R)
\]
by requiring that its fibers are given by $(M \otimes C)_{(b,t)}= M_{(b,t)} \otimes_{\IZ}  C_{(b,t)}$.  
In particular we consider the fiberwise complex $M \otimes F^{\delta}$.

Using  part of Remark~\ref{prove-hf} below and the fact that our fibers admit arbitrarily fine triangulations one can 
show that for $\delta \in \cals$ the complex $M \otimes F^{\delta}$ is homotopy equivalent inside 
$\ch \overline{ \calc }^{\Gamma} ( \tilde{E} \times \IT ; \cale )$ to a locally finite complex and hence
\[
M \otimes F^{\delta} \in  \ch_{\hf} \calc^{\Gamma} ( \tilde{E} \times \IT ; \cale ).
\]

\begin{remark} \label{prove-hf}
Let $T$ be a triangulation of the metric space $X=|T|$ such that the diameter of each simplex is smaller than $\delta$.
Let $C(T)$ denote the chain complex associated to the triangulation and let $C_{\sing}^{\epsilon} ( X )$ denote the 
subcomplex of the singular chain complex of $X$ generated by all singular simplices which are smaller than $\epsilon$. 
Both complexes can be considered as complexes over $X$ using the barycenters. If 
$0 < \delta \leq \delta_1 \leq \delta_2$ then there are natural inclusions
\[
C(T) \to C^{\delta_1} (X) \to C^{\delta_2} ( X ).
\]
Both maps are chain homotopy equivalences and one can show that the homotopy inverse and the homotopies can be chosen to
be $10 \delta_2$-controlled when considered as morphisms over $X$.
\end{remark}

\subsection{The transfer functor}

The discussion above suffices to define the desired transfer functor on objects. In order to define it on morphisms 
we need the fiber transport. For $(b,t)$ and $(b',t') \in \tilde{B} \times \IT$ the map
\[
\xymatrix{
E_{\pi(b)} & \ar[l]^-{\pi}_-{\cong} \tilde{E}_b \ar[r]^-{\nabla_{b',b}} & 
 \tilde{E}_{b'} \ar[r]_-{\pi}^-{\cong} & E_{\pi(b')}
         }
\]
induces a chain map
\[
F_{(b,t)} \to F_{(b',t')} 
\] 
which we will denote by $\nabla_{(b,t),(b',t')}$. Using 
Assumption~\ref{assumption-on-nabla}~\ref{assumption-functoriality} 
one checks that 
\begin{eqnarray*}
M & \mapsto  & M \otimes F \\
f=(f_{(b',t'),(b,t)}) & \mapsto &  f \otimes \nabla = ( f_{(b',t'),(b,t)} \otimes \nabla_{ (b',t'),(b,t)} )
\end{eqnarray*}
defines a functor to homotopy finite chain complexes over $\tilde{E} \times \IT$ if one ignores the control condition.
But since $\nabla_{b',b}$ can stretch simplices it does not induce a well defined map on the singular simplices of 
a fixed restricted size. This means that the analogous definition with $F^{\delta}$ does not work. In order to deal with this problem 
we formally enlarge our category. 

In Appendix~\ref{appendix-tilde-construction} we  construct 
for every Waldhausen category $\calw$ satisfying some mild conditions (which are satisfied
for categories of chain complexes in an additive category) a Waldhausen category 
$\widetilde{ \calw}$. Objects in this category are sequences
\[
\xymatrix{
C_0 \ar[r]^{c_0} &  C_1 \ar[r]^{c_1} &  C_2 \ar[r]^{c_2}  & \dots
         }
\]
where the $C_i$ are objects in $\calw$ and all $c_n$ are simultaneously cofibrations and weak equivalences 
(trivial cofibrations) in $\calw$. A morphism $f$ in $\widetilde{ \calw }$ is represented by a sequence 
$(f_m , f_{m+1} , f_{m+2} , \dots )$ of morphisms in $\calw$ which fit into a commutative diagram 
\[
\xymatrix{
C_m \ar[r]^{c_m} \ar[d]^{f_m} & 
C_{m+1} \ar[r]^{c_{m+1}} \ar[d]^{f_{m+1}} & C_{m+2} \ar[r]^{c_{m+2}}  \ar[d]^{f_{m+2}} & \dots \; \; \;  \\
D_{m + \kappa} \ar[r]^{d_{m + \kappa}}    & D_{m+ \kappa+1} \ar[r]^{d_{m+ \kappa+1}}     & 
D_{m+ \kappa +2} \ar[r]^{d_{m+ \kappa+2}} & \dots \; .
         }
\]
Here $m$ and $\kappa$ are nonnegative integers. If we enlarge $m$ or $\kappa$ the resulting diagram represents the 
same morphism, i.e.\ we identify  $(f_m , f_{m+1} , f_{m+2} , \dots )$ with the sequence 
$(f_{m+1} , f_{m+2} , f_{m +3 } , \dots )$ but also with 
$(d_m \circ f_m , d_{m+1} \circ f_{m+1} , d_{m+2} \circ f_{m+2} , \dots )$. Sending an object to the constant 
sequence yields an obvious inclusion $\calw \to \widetilde{\calw}$ and 
according to Proposition~\ref{enlarge-equivalence} 
this inclusion 
induces an equivalence in connective  $K$-theory. In the case where 
$\calw= \ch_{\hf} \calc^{\Gamma} ( \tilde{E} \times \IT ; \cale )$ we write 
$\widetilde{\ch}_{\hf} \calc^{\Gamma} ( \tilde{E} \times \IT ; \cale )$ for the corresponding enlargement and using
Remark~\ref{remark-only-connected} we conclude that the inclusion induces an equivalence in non-connective $K$-theory.

Let $\Jd = ( \delta^0 , \delta^1 , \delta^2 , \dots )$ be a monotone increasing sequence of decay speed functions 
$\delta^i \in \cals$, i.e.\ for all $i \geq 0$ and all $t \in \IT$ we have $\delta_t^i \leq \delta_t^{i+1}$.
Then for a module $M$ in $\calc^{\Gamma} ( \tilde{B} \times \IT ; \cale_d )$
\[
M \otimes F^{\Jd} = ( M \otimes F^{\delta^0} \to M \otimes F^{\delta^1} \to M \otimes F^{\delta^2} \to \dots )
\]
defines an object in $\widetilde{ \ch}_{\hf}  \calc^{\Gamma} ( \tilde{E} \times \IT ; \cale )$. Here the maps in the 
sequence are the natural inclusion maps. They are shown to be $\cale$-controlled homotopy equivalences using again 
Remark~\ref{prove-hf}.

\begin{assumption} \label{assumption-about-d}
Suppose for each $\alpha \geq 0$ there exists an integer $\kappa(\alpha)$ such that the following holds:
\begin{quote} 
If $d(b,b') \leq \alpha$ and  $| t - t'| \leq \alpha$ then for all $e \in \tilde{E}_{b'}$ we have 
\[
\nabla_{b,b'} ( \{ e \}^{\delta^i_t} ) \subset \nabla_{b,b'} ( \{ e \} )^{\delta^{i+\kappa( \alpha )}_{t'}},
\]
for all sufficiently large $t,t'$.
\end{quote}
Here the thickenings are taken in $\tilde{E}_b$ and $\tilde{E}_{b'}$ with respect to distance in
the ambient manifold.
\end{assumption}

Under this assumption we immediately obtain the following proposition.

\begin{proposition}
Suppose $\Jd$ is a sequence of decay speed functions satisfying Assumption~\ref{assumption-about-d}. 
Then there exists a functor 
\[
\tr^{\Jd}: \calc^{\Gamma} ( \tilde{B} \times \IT ; \cale_d ) \to 
        \widetilde{ \ch}_{\hf}  \calc^{\Gamma} ( \tilde{E} \times \IT ; \cale )
\]
which sends $f: M \to N$ to the morphism in
$\widetilde{ \ch}_{\hf}  \calc^{\Gamma} ( \tilde{E} \times \IT ; \cale )$ represented by 
\[
\xymatrix{
M \otimes F^{\delta_0} \ar[r] \ar[d]^-{f \otimes \nabla} & 
M \otimes F^{\delta_{1}} \ar[r] \ar[d]^-{f \otimes \nabla} & 
M \otimes F^{\delta_{2}} \ar[r] \ar[d]^-{f \otimes \nabla} & \dots \\
N \otimes F^{\delta_{\kappa}} \ar[r] & N \otimes F^{\delta_{\kappa +1}} \ar[r] & 
N \otimes F^{\delta_{\kappa +2}} \ar[r] & \dots
         }
\]
for suitably chosen $\kappa$  depending on the bound of $f$.
\end{proposition}

It remains to check, that in our situation where the fiber bundle is $S^+ \IH_{ \{ 0 \} } M  \to M$, $\nabla$ is the 
asymptotic fiber transport and $\cals=\cals_{\asy}$, we can find a suitable sequence $\Jd$ of decay speed functions.
The proof will use the fact that the fiber transport is H{\"o}lder continuous, compare Lemma~\ref{hoelder}.

\begin{lemma} \label{sequence-exists}
There exists a sequence $\Jd=( \delta^0 , \delta^1 , \delta^2 , \dots )$ of decay speed functions
$\delta^i \in \cals_{\asy}$ satisfying Assumption~\ref{assumption-about-d} with respect to the
asymptotic fiber transport.
\end{lemma}

\begin{proof}
We abbreviate the constant $\frac{a}{b}$ appearing in Lemma~\ref{hoelder} by $\lambda$.
Lemma~\ref{hoelder} implies that $\nabla_{b,b'}(e^\delta)$ is contained in 
$\nabla_{b,b'}(e)^{C_0(\alpha) \cdot \delta^{\lambda}}$. Thus our task is to find $\delta^i \in \cals_{asy}$ such
that for all $\alpha > 0$ there is $\kappa(\alpha)$ such that for all sufficiently large $t,t'$ with 
$|t-t'| \leq \alpha$ we have 
\[
C_0(\alpha) (\delta_t^i)^{\lambda} \leq \delta_{t'}^i.
\]
In the following it will be convenient to extend all functions on $\IT$ to functions on $\IR$ that are constant on 
$(-\infty,1]$. Start with an arbitrary $\delta^0 \in \cals_{asy}$. Choose inductively $\delta^i \in \cals_{asy}$ such 
that
\[
i \cdot (\delta^0_{t-i} + \dots + \delta^{i-1}_{t-i})^{\lambda} < \delta^i_t
\]
for all sufficiently large $t$. This is indeed possible because of Lemma~\ref{construct-S-asy}~\ref{construct-S-asy-2} and 
(A) and (B) before Definition~\ref{fol-control-w-decay}.
Choose now $\kappa(\alpha) \in \IN$ larger than $\alpha$ and $C_0(\alpha)$. Then
\[
C_0(\alpha) \cdot (\delta_t^i)^{\lambda} \leq (i + \kappa(\alpha)) \cdot (\delta_t^i)^\lambda 
    \leq \delta^{\kappa(\alpha)+i}_{t+i+\kappa(\alpha)}  
    \leq \delta^{\kappa(\alpha)+i}_{t'},  
\]
for sufficiently large $t,t'$ with $|t-t'|\leq \alpha$. (For the last inequality note that all functions in 
$\cals_{asy}$ are monotone decreasing.) 
\end{proof}

\subsection{An element in the Swan group} \label{subsection-mult-with-swan}

It remains to study the (non-commutative) triangle
\[
\xymatrix{
 & \widetilde{\ch}_{\hf} \calc^{\Gamma} ( \tilde{E} \times \IT , \cale )^{\infty} \ar[d]_-p \\
\calc^{\Gamma} ( \tilde{B} \times \IT , \cale_d )^{\infty}  \ar[ur]^-{\tr^{\Jd}}  \ar[r]^-{\inc}  & 
\widetilde{\ch}_{\hf} \calc^{\Gamma} ( \tilde{B} \times \IT , \cale_d )^{\infty} .
         }
\]
We will denote the induced maps in $K$-theory by the same symbols. Recall that the inclusion $\inc$ induces an 
isomorphism. What we would like to understand is the self map $\inc^{-1} \circ p \circ \tr^{\Jd}$. To describe the 
result we need some preparation. Let us fix a point $b_0 \in \tilde{B}$. Observe that in general the diagram
\[
\xymatrix{
\tilde{E}_{g b_0} \ar[dr]^-{\pi}_-{\cong} \ar[rr]^-{\nabla_{ b_{0} ,g b_{0}}} 
                             & & \tilde{E}_{b_{0}} \ar[dl]_-{\pi}^-{\cong} \\
& E_{\pi(b_0)} &
         }
\]
does not commute. In fact we use it to define a left $\Gamma$-operation on $E_{\pi(b_0)}$ by letting $g \in \Gamma$ 
act via $\pi \circ \nabla_{b_0 , gb_0} \circ \pi^{-1}$. The singular chain complex $F_0 = C_{\sing} ( E_{\pi ( b_0 )} )$
 hence becomes a complex of $\IZ \Gamma$-modules. As a $\IZ$-chain complex $F_0$ is homotopy equivalent to a finite 
complex of finitely generated $\IZ$-modules. The homology groups $H_i( F_0 )$ are hence $\IZ \Gamma$-modules which are 
finitely generated as $\IZ$-modules. Such a module defines an element in the Swan ring
$\Sw ( \Gamma ; \IZ )$. $K$-theory becomes a module over the Swan ring via maps
\[
\Sw ( \Gamma ; \IZ ) \otimes_{\IZ} K_n( R \Gamma ) \to K_n( R \Gamma ).
\]
The Swan ring, its action on $K$-theory and certain variants we need below in the proof are discussed in 
Appendix~\ref{appendix-swan}.

\begin{proposition} \label{mult-by-swan}
Under the identification 
\[
K_n ( R \Gamma ) \cong K_{n+1} ( \calc^{\Gamma} ( \tilde{B} \times \IT , \cale_d )^\infty)
\]
coming from the germs at infinity fibration (compare Example~\ref{germs-at-infty})
the map $\inc^{-1} \circ p \circ \tr^{\Jd}$ corresponds to multiplication with
\[
\Sigma_{i=0}^{\infty} (-1)^{i} \left[ H_i ( F_0 ) \right] \in \Sw ( \Gamma ; \IZ ).
\]
\end{proposition}

\begin{corollary} \label{corollary-triangle-commutes}
Diagram~\eqref{diagram-transfer} induces a commutative diagram in $K$-theory.
\end{corollary}

\begin{proof}
The fiber of $S^{+} \IH \tilde{M} \to \tilde{M}$ is contractible and hence the Swan group element is represented by the 
trivial $\IZ \Gamma$-module $\IZ$, which acts as the identity on $K$-theory.
\end{proof}

The rest of this subsection is devoted to the proof of Proposition~\ref{mult-by-swan}. Again we will only discuss the 
argument for connective $K$-theory. This yields Proposition~\ref{mult-by-swan} for $n \geq 1$. The general result 
follows by filling in extra $\IR^n$-factors, compare Remark~\ref{remark-only-connected}.

We first want to get rid of the $\widetilde{ \quad }$-construction. Consider
\[
\xymatrix{
\calc^{\Gamma} ( \tilde{B}\ \times \IT ,  \cale_d ) 
\ar@/^3ex/[rr]^-{p \circ \tr^{\Jd}}
\ar[r]_-{- \otimes \nabla}  & 
\ch_{\hf} \calc^{\Gamma} ( \tilde{B} \times \IT , \cale_d ) 
\ar[r]_-{\inc}  &
\widetilde{\ch}_{\hf} \calc^{\Gamma} ( \tilde{B} \times \IT , \cale_d ) 
         }.
\]
Here the functor $- \otimes \nabla$ is given by $M \mapsto M \otimes F$ and $f \mapsto f \otimes \nabla$, where now 
$F$ is considered as a complex over $\tilde{B} \times \IT$ and we do not care how it is distributed over each fiber.

\begin{lemma} \label{lemma-natural-transformation}
There is a natural transformation between $\inc \circ ( - \otimes \nabla )$  and $p \circ \tr^{\Jd}$ which is 
objectwise a weak equivalence.
\end{lemma}

\begin{proof}
At $M$ the natural transformation is given by the natural inclusion
\[
\xymatrix{
M \otimes F^{\delta_0} \ar[r] \ar[d] & M \otimes F^{\delta_1} \ar[r] \ar[d] & 
                                    M \otimes F^{\delta_2} \ar[r] \ar[d] & \dots \; \; \; \\
M \otimes F \ar[r]^-= &  M \otimes F \ar[r]^-= &  M \otimes F \ar[r]^-= & \dots \; .
         }
\]
\end{proof}

The functors and the natural transformation in Lemma~\ref{lemma-natural-transformation} are compatible with the germs 
at infinity fibrations. The middle terms in these fibrations are contractible since we work with the product metric.
(Compare Example~\ref{germs-at-infty} and \cite[Proposition~4.4, Example~4.5]{Bartels-Farrell-Jones-Reich(topology)} for such  arguments.)
Hence there is a version of Lemma~\ref{lemma-natural-transformation} for the germs at infinity categories. We have 
reduced our question to comparing the two maps
\[
\xymatrix{
\calc^{\Gamma} ( \tilde{B} ,  \cale_d ) 
\ar@<1ex>[r]^-{ \inc } 
\ar@<-1ex>[r]_-{- \otimes \nabla}  & 
\ch_{\hf} \calc^{\Gamma} ( \tilde{B} , \cale_d ) 
         }.
\]
Here $- \otimes \nabla$ is the obvious restriction of the functor above with the same name. Since we assume that 
$\tilde{B}$ is $\Gamma$-compact the $\cale_d$-condition is no extra condition and we omit it in the following. 
For the same reason the inclusion of the orbit $\Gamma b_0 \to \tilde{B}$ for some fixed $b_0 \in \tilde{B}$ induces 
equivalences on the categories and we are reduced to comparing the upper horizontal map in the following (non commutative) diagram to
the natural inclusion. (The diagram does commute if one replaces the horizontal maps by the natural inclusions.)
\begin{eqnarray} \label{diagram-square-up-to-natural}
\xymatrix{
\calc^{\Gamma} ( \Gamma b_0 ; R) 
\ar[r]^-{- \otimes \nabla}  \ar[d]_{ \oplus } & 
\ch_{\hf} \calc^{\Gamma} ( \Gamma b_0 ; R) \ar[d]_{ \oplus }  \\
\calc ( \pi (b_0 ) ; R \Gamma ) 
\ar[r]^-{- \otimes \id_{F_0} }  & 
\ch_{\hf} \calc ( \pi (  b_0 ) ; R \Gamma ).  
         }
\end{eqnarray}
The vertical functors in this diagram are equivalences given by sending an $R$-module $M=(M_{g b_0})$ over $\Gamma b_0$ 
to $\oplus_{g \in \Gamma} M_{g b_0}$ considered as an $R \Gamma$-module, 
compare Lemma~2.8 in~\cite{Bartels-Farrell-Jones-Reich(topology)}. 
The lower horizontal map sends an $R \Gamma$-module $N$ to the complex of $R \Gamma$ modules $N \otimes_{\IZ} F_0$. Here $F_0$ is 
the singular chain complex of the fiber $E_{\pi(b_0)}$ considered as a $\IZ \Gamma$-module as explained towards the 
beginning of this subsection and $\Gamma$ operates diagonally on $N \otimes_{\IZ} F_0$. On morphisms the functor sends 
$f$ to $f \otimes \id_{F_0}$. Diagram~(\ref{diagram-square-up-to-natural}) does {\em not} commute but we have the 
following lemma.

\begin{lemma} \label{another-nat-trafo}
There is a natural transformation between the two ways through Diagram~(\ref{diagram-square-up-to-natural})
which is objectwise an isomorphism.
\end{lemma}

\begin{proof}
Let $M = (M_{g b_0})$ be an $R$-module over $\Gamma b_0$. Both ways through 
(\ref{diagram-square-up-to-natural}) send $M$ to $(\oplus_{g \in \Gamma} M_{g b_0}) \otimes F$. However,
the $\Gamma$-actions are different. If we go first right and then down $\Gamma$ acts only on the first factor,
if we go down and then right $\Gamma$ acts diagonally on the tensor product. The natural transformation
from  right/down to down/right sends $m \otimes v \in M_{g b_0} \otimes F_0$ to 
$m \otimes g v \in M_{g b_0} \otimes F_0$.   
\end{proof}

It remains to compare $- \otimes \id_{F_0}$ to the natural inclusion.
In Appendix~\ref{appendix-swan} we explain a variant $\Sw^{\ch} (\Gamma ; \IZ )$ of the Swan group together
with its  action on $K$-theory. The complex $F_0$ is a complex of $\IZ \Gamma$-modules which is degreewise free as a $\IZ$-complex
and whose homology is finitely generated as an abelian group.
Such a complex
defines an element in $\Sw^{\ch} (\Gamma ; \IZ)$ and $\inc^{-1} \circ (- \otimes \id_{F_0})$
describes the action of this element on $K$-theory. Proposition~\ref{mult-by-swan} now follows from Proposition~\ref{all-three-isos}
and the discussion following that proposition.

\begin{remark} In the case we are interested in, where the fiber is a disk, 
we could avoid the Swan group and proceed differently after Lemma~\ref{another-nat-trafo}.
In that case the augmentation $\epsilon \colon F_0 \to \IZ$ is a $\Gamma$-equivariant 
map which is a non-equivariant chain homotopy equivalence.
In particular it induces for every free $R \Gamma$-module $N$ a homology isomorphism $N \otimes_{\IZ} F_0 \to N \otimes_{\IZ} \IZ$.
Each chain module of $N \otimes_{\IZ} F_0$ is non-canonically isomorphic to the $R \Gamma$-module with the same underlying abelian
group but where $\Gamma$ operates only on the left tensor factor. Hence both complexes are complexes of free $R \Gamma$-modules 
and the homology isomorphism is in fact an $R\Gamma$-chain homotopy equivalence. This yields a natural transformation between
$- \otimes F_0$ and the inclusion which is objectwise a weak equivalence.
\end{remark}


\typeout{---------------------fct -------------------------}
\section{A foliated control theorem for higher $K$-theory} \label{section-fct}

In this section we will show that a certain relax control map induces an equivalence in $K$-theory.
Roughly speaking the map relaxes control from metric control to foliated control with respect to the 
geodesic foliation. Unfortunately the precise statement is more complicated and to formulate it 
we need some rather lengthy preparations. 
The reader should right away take a look at Subsection~\ref{statement-fct} to get a first idea about the statement
we are aiming at.

Ignoring the technicalities the argument can be summarized as follows.
We know that metric control leads to homological behavior. In particular we have the long exact sequences associated to
pairs of spaces in order to work inductively over the skeleta of a cell structure. 
One task is now to formulate and prove  analogous results 
for foliated control using the skeleta of a long and thin cell structure.
The crucial step is a ``foliated excision'' result that reduces the statement about the relax
control map to a comparison result for a collection of long and thin cells 
(compare Proposition~\ref{foliated-excision}).
Carefully bookkeeping the error-terms
one can even assume that one has a collection of long and thin cells in Euclidean space
equipped with a ``standard'' $1$-dimensional foliation. The comparison result is then easily established:
an Eilenberg swindle is used to reduce the question to transversal cells. On transversal cells
metric control and foliated control coincide (compare Lemma~\ref{euclidean-standard}).


\subsection{Flow cell structures} \label{subsection-flow-cell-structures}

Let $N$ be an $n$-dimensional  Riemannian manifold which is equipped with a 
smooth flow $\Phi$. 
The flow determines a one-dimensional foliation which will be called $F$.
In our application $N$ will always be $S \IH {M}$ (or its universal covering) 
equipped with the 
geodesic flow and the corresponding foliation $F_{geo}$.

The following definition of a flow cell combines Definition~7.1 and Lemma~8.1 in \cite{Farrell-Jones(dynamics-I)}. Whereas in 
\cite{Farrell-Jones(dynamics-I)} Lemma~8.1 the information about the length of a cell is contained in the map $g_e$
we require the map to roughly preserve the length 
and use  instead a long  parametrizing interval $A_e$.
Below $\IR^n = \IR \times \IR^{n-1}$ is equipped with the standard Euclidean metric and foliated by the lines parallel
to the first coordinate axis. We denote this foliation by $F_{\IR^{n}}$. 
Moreover $\mu_n =10^{n+3}$ is the constant which depends only on the dimension 
that appears in Proposition~7.2 in \cite{Farrell-Jones(dynamics-I)}.
Recall that for a subset $Y\subset \IR^n$ we denote by $Y^{(\alpha , \delta)}$ the set of all $x \in \IR^n$ for 
which there is a $y \in Y$ such 
that $d_{F_{\IR^{n}}}(x,y) \leq (\alpha,\e)$, compare Subsection~\ref{section-fol-control-w-decay}.

\begin{definition}[$\beta$-flow cell] \label{definition-flow-cell}
Let $\beta>0$ be given.
A cell $e \subset N$ is called a {\em $\beta$-flow cell} if there exist cells 
$A_e \subset \IR$, $B_e \subset \IR^{n-1}$, a number  $\epsilon_e >0$ and a smooth
embedding $g_e\colon ( A_e \times B_e)^{(\beta,\epsilon_e)} \to N$ such that
\begin{enumerate}
\item
We have $g_e( A_e \times B_e )=e$.
\item
The map $g_e$ preserves the foliation, i.e.\ for each $y \in \IR^{n-1}$
the segment $\IR \times \{ y \} \cap (A_e \times B_e)^{(\beta , \epsilon_e)}$
is mapped to a segment of a flow line in $N$.
\item
If $A_e \subset \IR$ is not a $0$-cell then it is an interval of length exactly $\beta$.
\item \label{dfc-vier}
For all tangent vectors $v \in T ((A_e \times B_e)^{( \beta , \epsilon_e)} )$ which are tangential to
the flow lines we have
\[
| v | < | dg_e ( v ) | \leq \mu_n / 5 \cdot |v|.
\]
\end{enumerate}
\end{definition}

There are two sorts of flow cells: 
A flow cell where $A_e$ is a $0$-cell will be called {\em transversal}.  
If $A_e$ is a $1$-cell we call the cell $e$ a {\em long} cell.
Observe that from (iii) and (iv) it follows that such a cell 
is $\beta$-long in the sense that for every $y \in B_e$ the 
segment $g_e ( A_e \times \{ y \} )$ has arclength strictly 
larger than $\beta$ (and shorter than $\beta \cdot \mu_n/5$).

\begin{remark}  \label{remark-flow-cell}
A simple compactness argument shows that for a flow cell we additionally have the following.
\begin{enumerate}
\item[(v)]
There exists a constant $C_e>1$ such that 
for all tangent vectors $v \in T((A_e \times B_e)^{( \beta , \epsilon_e)} )$ we have 
\[
C_e^{-1} \cdot | v | \leq | dg_e ( v ) | \leq C_e \cdot | v |.
\]
\end{enumerate}
\end{remark}

\begin{remark} \label{compare-foliated-distance}
Since (iv) and (v) hold over the $(\beta , \epsilon_e)$-thickening of $A_e \times B_e$,
foliated distances (compare Subsection~\ref{section-fol-control-w-decay}) between 
points in the cell which are small compared to $(\beta , \epsilon_e)$
can be approximately determined in Euclidean space using the chart. More precisely:
Given $z=g_e(v)$ and $z'=g_e(v')$ with $z,z' \in e$ and  
\[
d_F ( z , z') \leq (\alpha , \delta) \leq ( 5 \mu_n^{-1} \cdot \beta , C_e^{-1} \cdot \epsilon_e )
\]
we have 
\[
d_{F_{\IR^n}} ( v , v' ) \leq ( \alpha , C_e \delta ). 
\]
The other way round 
\[
d_{F_{\IR^n}} ( v , v' ) \leq ( \alpha , \delta )  \quad \mbox{ implies } \quad
d_F ( z , z') \leq ( \mu_n \cdot \alpha , C_e \cdot \delta) .
\]
\end{remark}

A cell structure $L$ for a compact subset of $N$ all whose cells are $\beta$-flow cells 
will be called a $\beta$-\emph{flow cell structure}. 
Given a $\beta$-flow cell structure we will always fix choices of charts $g_e$ and constants $\epsilon_e$ and $C_e$ as in 
Definition~\ref{definition-flow-cell}.
For a given cell structure $L$ we denote by $|L| \subset N$ its underlying topological space.
We recall the main result of Proposition~7.2 and Lemma~8.1 in \cite{Farrell-Jones(dynamics-I)}.
 
\begin{theorem} \label{existence-cell-structure}
Let $n=\dim N$ and $\mu_n=10^{n+3}$. Let $N^{\leq \mu_n \beta}$ denote the union of all leaves which are 
shorter than $\mu_n \beta$.
For arbitrarily large $\beta$ and any 
compact subset  $K \subset N -  N^{ \leq \mu_n \beta }$ there exists a $\beta$-flow-cell structure $L$ with $K \subset |L|$.
\end{theorem} 

Given a cell structure $L$ we denote by $L^{[k]}$ the set of all $k$-cells and  by $L^{(k)}$ the set of 
cells of dimension less than or equal to $k$. The $k$-skeleton is $|L^{(k)}|$.
We define a filtration  
\[
N^{(-1)} \subset N^{(0)} \subset \dots \subset N^{(n)} =N
\]
of $N$ as follows. Set $N^{ (-1) }= \overline{ N - |L| }$
and $N^{ (k) }= N^{ (-1) } \cup | L^{(k)} |$. Observe that a cell $e \in L^{ [k] }$ may already be contained
in $N^{ (-1) }$ and hence does not contribute to the $k$-th filtration step. We hence define $L^{ \{ k \} }$
to be the set of those $k$-cells which do not lie entirely in $N^{ (-1) }$. Note that such a cell $e$
can meet $N^{(-1)}$ only with its boundary, which we denote $\partial e$. 

Let $\tilde{N}$ denote the universal cover of $N$ and let $\Gamma$ be the fundamental group
which acts via deck transformations on $\tilde{N}$. 
The lifted cell structure will be denoted $\tilde{L}$ and 
$\tilde{N}^{(k)}$ denotes the preimage filtration of $\tilde{N}$ under the covering projection.
Also we will use $\tilde{L}^{ [k] }$ and  $\tilde{L}^{ \{ k \} }$ to denote the obvious sets of cells of $\tilde{L}$.
Note that a flow cell in $N$ gives rise to a whole $\Gamma$-orbit of flow cells in $\tilde{N}$ for 
which one can simultaneously
use the constants appearing in Definition~\ref{definition-flow-cell} and Remark~\ref{remark-flow-cell}.

Roughly speaking the following proposition says that if one removes an $(\alpha , \Delta )$-thickening 
of the $(k-1)$-st filtration step from the $k$-th filtration step, then the remaining
pieces of the interiors of the $k$-cells are at least $(\alpha , \delta)$-foliated apart
from one another. This fact will later play a crucial role in Proposition~\ref{foliated-excision}.

\begin{proposition} \label{flow-cell-structure-constants}
Let $L$ be a $\beta$-flow cell structure for a compact subset $|L| \subset N$. 
There exists an $\epsilon_L^{\prime} >0$ and 
a function $\Delta_L( \alpha , \delta )$
defined for $0 \leq \alpha \leq \beta$ 
and $0 < \delta < \epsilon_L^{\prime}$
such that $\Delta_L( \alpha , \delta ) \geq \delta$ and the following holds:
\begin{enumerate}
\item \label{fcsc-eins}
Suppose  $e \in \tilde{L}^{ \{ k \} }$ is a $k$-cell and $x \in \tilde{N}^{ (k) }$. 
Whenever $\Delta_L ( \alpha , \delta )$ is defined and 
$\Delta \geq \Delta_L ( \alpha , \delta )$ then
\[
d_F( e - \partial e^{(\alpha,\Delta)} , x  ) \leq ( \alpha , \delta )
\]
implies that $x \in e$.
\item \label{fcsc-zwei}
For fixed $\alpha$ the function $\Delta( \alpha , \delta )$ 
tends monotone to zero when $\delta$ does.
\end{enumerate}
\end{proposition}

\begin{remark}
Observe that Proposition~\ref{flow-cell-structure-constants} says in particular that for 
two $0$-cells $e$ and $e^{\prime} \in \tilde{L}^{ \{ 0 \} }$ with
\[
d_F ( e , e^{\prime} ) < ( \beta , \epsilon_L^{\prime} )
\] 
we have $e = e^{\prime}$ since $\partial e = \emptyset$ for $e \in \tilde{L}^{[0]}$.
\end{remark}

\begin{proof}[Proof of Proposition~\ref{flow-cell-structure-constants}] 
For $e \in \tilde{L}^{ \{ k \} }$, $\Delta > 0$ and 
$\alpha \geq 0 $ let $Y(e,\alpha,\Delta) \subset \tilde{N}$
consist of all points $x = g_e(t_0+t,y)$ where 
$g_e(t_0,y) \in  e - \partial e^{(\alpha,\Delta) }$ and 
the path $\tau \mapsto g_e(t_0 + \tau,y)$ for $0 \leq \tau \leq t$ (resp. $t \leq \tau  \leq 0$) 
has arclength  $\leq \alpha$. 
(If the flow has unit speed then 
$Y(e, \alpha , \Delta )$ coincides with 
$\Phi_{[- \alpha , \alpha]} ( e - \partial e^{(\alpha,\Delta) } )$.)
If $\alpha \leq \beta$ then $Y(e,\alpha,\Delta)$ is disjoint from
$\tilde{N}^{(k-1)}$ 
and every cell $e' \in \tilde{L}^{ [k] }$ 
unless $e=e'$. For a long cell $e$ this is immediate from the construction (and the fact that we assume $\Delta >0$). 
To see it in the case where
$e$ is a transversal cell observe that a cell in $\tilde{L}^{ \{ k \} }$ can meet $\tilde{N}^{ (-1) }$
only at its boundary and hence points near (but not in) $e - \partial
e^{(\alpha , \Delta)}$ 
which lie on a flow line which  meets
$e - \partial e^{ ( \alpha , \Delta) }$ must lie
in the interior of a long cell of dimension bigger than $k$. (Here one
uses the fact that by definition the 
long cells are \emph{strictly}
longer than $\beta$.)
Define $X ( e )= \tilde{N}^{(k-1)} \cup \{ |e^{\prime}| \; | \; e^{\prime} \in L^{ \{ k \} } , \; e^{\prime} \neq e \}$.
As a first approximation to the foliated distance appearing in (i) we discuss the distance
\[
d(e, \alpha , \Delta ) = d ( Y ( e , \alpha , \Delta ) , X( e ) ).
\]
Observe that even though $X(e)$ is usually not compact only the intersection of it with some sufficiently 
large compact set matters. For small enough $\Delta$ we know that $Y ( e , \alpha , \Delta ) \neq \emptyset$
and hence $d( e  , \alpha , \Delta )$ is a positive number, say $4 \epsilon^{\prime}_e$.
Moreover $\Delta \leq \Delta^{\prime}$ implies $d(e, \alpha , \Delta) \leq d(e , \alpha , \Delta)$ and
$d( e , \alpha , \Delta)$ tends to $0$ if $\Delta$ does unless $e$ is a $0$-cell.
For cells $e$ which are not $0$-cells and $0 < \delta <
\epsilon_e^{\prime}$ we define 
$\Delta_e ( \alpha , \delta)$ as the minimal $\Delta$
for which $d(e, \alpha , \Delta ) \geq 2 \delta$. 
Since each deck transformation $\gamma \in \Gamma$ acts by isometries, preserves the foliation, respects the filtration
and permutes the cells, we know that  $d( e, \alpha , \Delta ) = d ( \gamma e , \alpha , \Delta )$.
Since there are only finitely many orbits of cells we can define
$\epsilon^{\prime}$ 
as the minimal $\epsilon_e^{\prime}$, where
$e$ ranges over all cells which are in $L^{\{ k \} }$ for some $k \geq 0$
and $\Delta ( \alpha , \delta )$ as the maximal $\Delta_e ( \alpha , \delta )$, where $e$ ranges over all cells which are 
in $L^{ \{ k \} }$ for some $k \geq 1$.
For $0 \leq \alpha \leq \beta$ and $0 < \delta < \epsilon_L^{\prime}$ and every $e \in L^{ \{ k \} }$ 
we have $d(e, \alpha , \Delta ) \geq 2 \delta$ for all $\Delta \geq \Delta ( \alpha , \delta )$ and 
$\Delta ( \alpha , \delta )$ tends to $0$ if $\delta$ does. 
It remains to improve the established inequalities slightly.
Note first that for $\alpha \geq 0$ there is a constant $C_{flw}(\alpha)$ 
such that for $x$, $y \in |\tilde{L}|$ with $d(x,y) \leq 1$ and all $t$ with $-\alpha \leq t \leq \alpha$ we have
$d(\Phi_t(x),\Phi_t(y)) < C_{flw}(\alpha) \cdot d(x,y)$ by a compactness
argument. (For the geodesic flow on $S \IH \tilde{M}$
this holds even over all of $S \IH \tilde{M}$ by \ref{flowconstant}.) 
Because of our symmetric 
definition of foliated distance in Subsection~\ref{section-fol-control-w-decay} we see that
$\Delta_L(\alpha,\delta) = \Delta(\alpha,(C_{flw}(\alpha)+2) \cdot \frac{\delta}{2} )$ 
and
$\e'_L = \min(1,(C_{flw}(\beta)+1)^{-1} \cdot \frac{\e'}{2})$ satisfy our requirements.
Compare also Lemma~\ref{foltriangle} \ref{foltriangle-eins}.
\end{proof}


\subsection{A family of flow cell structures} \label{subsection-famil-of-flow}

In order to prove foliated control results for $S \IH \tilde{M}$ equipped with the geodesic 
foliation $F_{geo}$ we need 
longer and longer cell structures (necessarily missing more and more closed  geodesics) 
but we also want to cover larger and larger chunks in the non compact $\IH$-direction 
(because the flow moves things in that direction).
This naturally leads us to choose flow-cell structures $L_{\beta , i}$ indexed by 
$\IN_0 \times \IN$  which are $\beta$-long
and cover the $[ -i , i ]$-part in the $\IH$-direction (each individual cell structure will 
only cover a compact region). Here are the details:

Let $\beta>0$ be given. 
Let $S \IH M^{\leq \mu_n \beta}$ denote the subset of $S \IH M$ that consists of all closed 
geodesics of length $\leq \mu_n \beta$. Here $\mu_n= 10^{n+3}$ with 
$n=\dim S \IH M$, cf.~Theorem~\ref{existence-cell-structure}.
Observe that $S \IH M^{\leq \mu_n \beta}$ lies in $S \IH_{ \{ 0 \} }
M$ 
because all compact flow lines have an $\IH$-coordinate
which is constantly $0$.

For a fixed $\beta \in \IN_0$ we choose a monotone decreasing 
sequence of tubular 
neighborhoods $T_{\beta , i }$, $i \in \IN$ of 
$S \IH M^{\leq \mu_n \beta }$ such that 
\[
\bigcap_{i \in \IN} T_{\beta , i } = S \IH M^{ \leq \mu_n \beta}.
\]
We will use the tilde-notation, i.e.\ $S\IH \tilde{M}^{\leq \mu_n \beta}$ and $\tilde{T}_{\beta , i}$, to
denote the obvious preimages under the universal covering projection $S\IH \tilde{M} \to S \IH M$.
Throughout the rest of this section we also 
fix a choice of a $\beta$-flow cell structure $L_{\beta, i}$ for $\beta \in \IN_0$ and $i \in \IN$
such that 
\[
S \IH_{[-i , i]}  M - T_{\beta,i} \subset | L_{\beta , i} | \subset S \IH  M,
\]
compare Theorem~\ref{existence-cell-structure}. 
Since its cells are shorter than $\beta \cdot \mu_n /5$ we can also arrange that
\[
| L_{\beta , i} | \subset S \IH _{[-i- \mu_n \beta , i + \mu_n \beta ]} M.
\]

To our choice of flow cell structures $L_{\beta , i}$  we will now associate certain sequences of constants and functions.
First recall that for a cell $e$ in a single flow cell structure $L$ we have the 
constants $\epsilon_e$ and $C_e$ appearing in Definition~\ref{definition-flow-cell}
and Remark~\ref{remark-flow-cell}. Moreover there are  the constant $\epsilon_L^{\prime}$ and the function 
$\Delta_L ( \alpha , \delta )$
from Proposition~\ref{flow-cell-structure-constants}.
We set
\begin{eqnarray} \label{constants-for-single-L}
\epsilon_L  =  \min \{ \epsilon_e | e \in L \} \cup \{ \epsilon_L^{\prime}  \} \quad 
\mbox{ and } \quad C_L  =  \max \{ C_e | e \in L \}.
\end{eqnarray}

Now back to the family 
$(L_{\beta , i})_{(\beta , i) \in \IN_0 \times \IN}$ we have chosen above. We set
\begin{eqnarray}
\e_i & = & \min \{ \e_{L_{\beta,i}} | 0 \leq \beta \leq i \},  \label{epsiloni} \\
C_i & = & \max \{ C_{L_{\beta,i}}  | 0 \leq \beta \leq i \}, \label{C-i} \\
\Delta_i(\delta) & = 
  & \max \{ \Delta_{L_{\beta,i}} (\alpha , \delta) | \alpha \in \IN_0 \mbox{ and } 0 \leq \alpha \leq \beta \leq i \}  
\mbox{ for } \delta < \e_i.  \label{alpha-in-N}
\end{eqnarray}   
Note that $\Delta_i( \delta ) \geq \delta$ and for $i$ fixed $\Delta_i ( \delta )$ tends to $0$ with $\delta$, compare 
Proposition~\ref{flow-cell-structure-constants}~\ref{fcsc-zwei}.  For $\delta$ fixed
$\Delta_i ( \delta )$ is monotone increasing with $i$.
Making the $\epsilon_i$ smaller and the $C_i$ bigger if necessary we will assume that $\epsilon_i$ tends monotone to $0$
and the $C_i$ form an increasing sequence of numbers $>1$.

Later on, we will be in a situation where we can ignore all $L_{\beta,i}$ with $\beta >  i$.
(Compare Proposition~\ref{first-reduction} and the definition of $\calf_{\IT( \beta )}$ before that proposition.)
With the above definitions the constants and functions labeled with $i$ have the desired properties simultaneously 
for all cell structures $L_{\beta, i}$ with $\beta \leq i$.

More precisely we have the following lemma. 
\begin{lemma} \label{simultaneous-for-all} \hfill
\begin{enumerate}
\item \label{simultaneous-for-all-C}
For fixed $i$
Remark~\ref{compare-foliated-distance} with $\epsilon_i$ and $C_i$ instead of $\epsilon_e$ and
$C_e$ applies simultaneously to all cells in all the cell structures $L_{\beta , i }$ 
with $\beta \leq i$.
\item \label{simultaneous-for-all-Delta}
Similarly for a fixed $\alpha \in \IN_0$ Proposition~\ref{flow-cell-structure-constants}~(i) 
applies with $\epsilon_i$ instead of $\epsilon_L^{\prime}$ to all cell
structures $L_{\beta , i}$ 
with $\alpha \leq \beta \leq i$.
\end{enumerate}
\end{lemma}


\subsection{Construction of the decay speed $\cals$} \label{subsection-construction-of-S}

Let $\Jt = (t_1 ,t_2 , \dots )$ be a sequence of numbers with $t_1=1$ and 
$t_i < t_{i+1}$. Given a sequence $( \delta_i )$ we define the associated 
step-function $\step_{\Jt} ( (\delta_i))$ to be the function on $\IT$ whose value on 
the interval $\left[ t_i , t_{i+1} \right)$ is $\delta_i$.
This defines a map from the space of sequences to the space of functions.

Our aim is now to construct a certain set of sequences $\calt= \{ (\delta_i) \}$ 
which will then
(after a choice of a suitable sequence $\Jt$) 
lead to the set of functions $\cals=\step_{\Jt} ( \calt )$ 
used to describe the decay speed in the $\cale_w$-control condition. 
In fact because 
of Remark~\ref{only-need-germs}
we are really only interested in the germs at infinity of such sequences.

\begin{lemma} \label{construct-T}
There exists a nonempty set $\calt=\{ ( \delta_i)_{i \in \IN} \}$ of sequences of 
positive numbers (each of which tends to zero) satisfying: 
\begin{enumerate}
\item \label{cT-vier}
Each sequence $( \delta_i) \in \calt$ is eventually smaller than the sequence 
$(\epsilon_i)$ defined in equation (\ref{epsiloni}), i.e.\ there exists an $i_0 \in \IN$ such that 
$\delta_i < \epsilon_i$ for all $i \geq i_0$.
\item \label{cT-fuenf}
For every $(\delta_i) \in \calt$ the sequence $(\Delta_i ( \delta_i ) )$,  
eventually defined by (i), lies again in $\calt$.
\item \label{cT-drei}
For every $(\delta_i) \in \calt$ the sequence $(C_i \cdot \delta_i)$ lies again in $\calt$.
\end{enumerate}
Moreover we have the following more elementary properties corresponding 
to (A) and (B) before Definition~\ref{fol-control-w-decay}.
\begin{enumerate}
\item[(A)] \label{cT-eins}
For $(\delta_i) \in \calt$ and $k \in \IZ$ we have 
$(\delta_{i+k})_{i \in \IN} \in \calt$. (Here we set $\delta_{i+k} = \delta_1$ for $i+k \leq 0$.)
\item[(B)] \label{cT-zwei}
Given $(\delta_i)$, $(\delta_i^{\prime}) \in \calt$ there exists 
$(\delta_i^{\prime \prime}) \in \calt$ such that 
$\delta_i + \delta_i^{\prime} \leq \delta_i^{\prime \prime }$ for all $i \in \IN$.
\end{enumerate}
\end{lemma}

\begin{proof} 
It will be convenient to define
$\Delta'_i(\delta)=\max \{ 2 \delta, \Delta_i(\delta), C_i \cdot \delta \}$
for $\delta \leq \e_i$ and $\infty$
otherwise.
For fixed $i$ this tends monotone to $0$ with $\delta$, for fixed $\delta$ it is monotone increasing with $i$.
The space $\calt$ consisting of all $(\delta_i)$ satisfying the
following condition
\[
\forall k,l \in \IZ,j \in \IN_0 \hspace{3ex} \exists i_0 \hspace{3ex}
\forall i \geq i_0
      \hspace{3ex} (\Delta'_{i+l})^{\circ j}(\delta_{i+k}) < \e_i
\]
is nonempty and satisfies (i), (ii), (iii), (A) and (B):
Properties (i) and (A)
are clear from this construction.
To check property (ii) observe that
\[
(\Delta'_{i+l})^{\circ j} (\Delta_{i+k}(\delta_{i+k}))
\leq
\left\{
  \begin{array}{rcl}
    (\Delta'_{i+l})^{\circ j+1} (\delta_{i+k}) & \mbox{ if } & i+l \geq i+k \\
    (\Delta'_{i+k})^{\circ j+1} (\delta_{i+k}) & \mbox{ if } & i+l \leq i+k
  \end{array}
\right.   .
\]
Property (iii) follows by replacing $\Delta'_{i+k}$ by $C_{i+k}$ in this inequality.
To check (B) first observe that
we have $(\Delta'_i)^{\circ j}(2 \delta) \leq (\Delta'_i)^{\circ
j+1}(\delta)$.
Thus,
$(\delta_i) \in \calt$ implies $(2 \delta_i) \in \calt$.
Moreover $(\delta_i) \in \calt$ and $\delta'_i \leq \delta_i$ implies
$\delta'_i \in \calt$. Finally, by construction,
$(\delta_i),(\delta'_i) \in \calt$ implies $\max(\delta_i,\delta'_i) \in
\calt$.
So we get (B) since
$\delta_i + \delta'_i \leq 2 \max(\delta_i,\delta'_i)$.
Since $\Delta^{\prime}_i( \delta )$ tends to $0$ with $\delta$ we can find
$\delta_i$ such that
$(\Delta'_{2i})^{\circ 2i}(\delta_i) < \e_{2i}$.
Then
\[
(\Delta'_{i+l})^{\circ j} (\delta_{i+k}) \leq (\Delta'_{i+l})^{\circ i+l} (\delta_{i+k}) 
 \leq (\Delta'_{2(i+k)})^{\circ 2(i+k)}(\delta_{i+k}) \leq e_{2(i+k)}   < \e_i,
\]
for sufficiently large $i$, i.e. if $i+l > j$, $1 \leq i+l \leq 2(i+k)$
and $i \leq 2(i+k)$.
Thus $\calt$ contains $(\delta_i)$ and is indeed not empty.
\end{proof}

Now choose an increasing sequence $\Jt=(t_1, t_2 , \dots )$ with $t_1=1$ such that
$t_{i+1} - t_i \geq i$ (this will be important in Proposition~\ref{second-reduction})
and $(e^{-at_i}) \in \calt$. Here $a$ is the curvature bound from (\ref{curvature-constants}).
The next statement is immediate from Lemma~\ref{construct-T}. 
 
\begin{proposition} \label{S-properties}
The set $\cals=\step_{\Jt} \calt$ satisfies the 
standard properties (A) and (B) of a class of decay speed functions introduced before
Definition~\ref{fol-control-w-decay}, each $\delta_t \in \cals$ tends to zero for $t \to \infty$ and moreover:
\begin{enumerate}
\item \label{S-properties-eins}
For each $A \cdot e^{-at} \in \cals_{geo}$ there exists a 
$\delta_t \in \cals$ with $A \cdot e^{-at} \leq \delta_t$, i.e.\
we can ``relax control'' from  $\cals_{geo}$ to  $\cals$.
\item
For a given $\delta_t=\step_{\Jt} (\delta_i) \in \cals=\step_{\Jt} \calt$ 
take $i_0 \in \IN$ as in Lemma~\ref{construct-T}(i)
and define for $t \geq t_{i_0}$ the function 
$\Delta_t = \step_{\Jt}  ( \Delta_i ( \delta_i ))$. This function
lies again in $\cals$.
\item \label{S-properties-drei}
For a given $\delta_t = \step_{\Jt} ( \delta_i) \in \cals $ the function
$C_t \cdot \delta_t = \step_{\Jt} ( C_i \cdot  \delta_i )$ lies again in $\cals$.
\end{enumerate}
\end{proposition}

Condition (i) is important to obtain the map (5) in our main diagram in Section~\ref{section-proof}, i.e.\
to connect up the following constructions with
the kind of control we obtained via the geodesic flow. 
The second and third condition will play an important role in
Proposition~\ref{foliated-excision} respectively in Lemma~\ref{compare-euclidean} below.


\subsection{Statement of the Foliated Control Theorem} \label{statement-fct}

Now we are prepared to define the control structures $\cale_w$ and $\cale_s$ on $S \IH \tilde{M} \times \IB \times \IT$
which were already mentioned in
the outline of the proof in Section~\ref{section-proof}.
Let us recall what happened so far. 
For all natural numbers 
$\beta \geq 0$ and $i \geq 1$ we chose
a tubular neighborhood $T_{\beta , i}$ of $S \IH \tilde{M}^{\leq \mu_n \beta}$, and  
$\beta$-flow cell structures $L_{\beta , i}$ such that 
$S \IH_{[-i,i]} M - T_{\beta,i} \subset |L_{\beta , i}|$.
The associated constants $C_{L_{\beta,i}}$, the $\epsilon_{L_{\beta , i}}$ 
(see (\ref{constants-for-single-L}), Remarks~\ref{remark-flow-cell}
and \ref{compare-foliated-distance})
and the functions $\Delta_{L_{\beta , i}}( \alpha, \delta )$ given by
Proposition~\ref{flow-cell-structure-constants} 
were used to define the 
sequences $\epsilon_i$, $C_i$ and $\Delta_i( \delta )$.
In Lemma~\ref{construct-T} we produced a space of sequences $\calt$ out of this data. 
Before Proposition~\ref{S-properties}
we then chose a sequence $\Jt=(t_1,t_2,\dots )$ and defined the set of decay speed functions $\cals= \step_{\Jt} ( \calt )$. 

Now set 
\begin{eqnarray} \label{definition-S}
S = \bigcup_{\beta , i} |L_{\beta , i}| \times [ \beta , \beta+1 ] \times [ t_i, t_{i+1} ] 
         \subset S \IH \tilde{M} \times \IB \times \IT.
\end{eqnarray}
The Foliated Control Theorem~\ref{fct} will improve control precisely over $S$.
Equip $\IB$ with the Euclidean metric and $S \IH \tilde{M} \times \IB$ with some product metric. 
Equip $\IB$ with the $0$-dimensional foliation by points and let $F_w$ denote the product 
foliation with the foliation $F_{geo}$ on $S \IH \tilde{M}$ given by the geodesic flow. 

\begin{definition}[Weak and strong control] \label{definition-weak-strong-control} \hfill
\begin{enumerate}
\item
The ``weak'' morphism control condition $\cale_w$ on $S \IH \tilde{M} \times \IB \times \IT$ is 
defined as foliated control with respect to the foliation $F_w$ on 
$S \IH \tilde{M} \times \IB$  with decay speed $\cals= \step_{\Jt}( \calt )$, 
i.e.\ $\cale_w = \cale( S \IH \tilde{M} \times \IB , F_w , \cals )$, cf.~Definition~\ref{fol-control-w-decay}.
\item
The ``very strong'' morphism control condition $\cale_{vs}$ on $S \IH \tilde{M} \times \IB \times \IT$
is defined as metric control with respect to the product metric on $S \IH \tilde{M} \times \IB$ with
decay speed $\cals= \step_{\Jt} ( \calt )$, i.e.\ $\cale_{vs} = \cale(S \IH \tilde{M} \times \IB , \cals)$, 
cf.~Definition~\ref{fol-control-w-decay}.
We let $\cale_{vs}^{\prime}$ denote the metric $\cale_{vs}$-control
over  the  subset $S \subset S \IH \tilde{M} \times \IB \times \IT$. 
(This is a control condition over $S \IH \tilde{M} \times \IB \times \IT$ as
explained in Definition~\ref{control-over}.)
\item
The ``strong'' morphism control condition $\cale_s$ is defined as $\cale_w \cap \cale_{vs}^{\prime}$,
i.e. foliated control everywhere and metric control over $S$.
\end{enumerate}
\end{definition}

Observe that strong and weak control differ only over the subset $S$, 
where $\cale_s$ requires the stronger metric control instead of only foliated control.

\begin{remark}
For the foliation $F_w$ there is a foliated triangle inequality analogous to 
Lemma~\ref{foltriangle}(i). Together with 
(A) and (B) (compare Proposition~\ref{S-properties}) it hence follows that $\cale_w$ is closed 
under composition and indeed defines a morphism control condition, compare Warning~\ref{warning-fol-control}.
\end{remark}

Recall that in  Section~\ref{section-proof} we introduced the
object support condition 
\begin{eqnarray*}
\calf_{\IB}    & = & \{ \; \{ (v, \beta , t) \; | \; \beta \leq \beta_0 \} \; | \; \beta_0 \in \IB \}  
\end{eqnarray*}
on $S \IH \tilde{M} \times \IB \times \IT$ 
and the subspace $(S \IH \tilde{M} \times \IB \times \IT)_\angle$ consisting of all
$(v,\beta,t)$ with $|h(v)| \leq t + \mu_n \beta$.
Here $h$, $\beta$ and $t$ denote the $\IH$-, $\IB$- respectively $\IT$-coordinate of a point 
$(v, \beta, t) \in S \IH \tilde{M} \times \IB \times \IT$ and $\mu_n=10^{n+3}$ with $n=\dim S \IH \tilde{M}$.
After all these preparations we can finally formulate the main result of this Section.

\begin{theorem}[Foliated Control Theorem] \label{fct}
The forget control map
\[
\calc^{\Gamma} ( (S \IH \tilde{M} \times \IB \times \IT)_\angle , \cale_s , 
                                  \calf_{\IB} )^{\infty}
\to 
\calc^{\Gamma} ( (S \IH \tilde{M} \times \IB \times \IT)_\angle , \cale_w , 
                                  \calf_{\IB} )^{\infty}
\]
given by relaxing the $\cale_s$-control condition to the $\cale_w$-condition induces an 
equivalence in $K$-theory.
\end{theorem}

After two preliminary reduction steps this result will be proven by induction over the 
skeleta of a relative cell structure. The proof will occupy the rest of this section.


\subsection{First reduction -- Delooping in the $\IB$-direction} 
\label{subsection-first-reduction}

Suppose we are given an $(\alpha , \delta_t)$-controlled morphism. Then in a region of the space $S \IH \tilde{M} \times \IB \times \IT$
where $\beta$ is larger than $\alpha$ and
$t$ is very large the morphism is quite well adapted to the flow cells we will find there. Conversely we would like to ignore
a certain region where we have no hope to prove a comparison result between foliated and metric control.
To capture this idea we introduce further object support conditions. We define 
analogous to $\calf_{\IB}$ the following object support conditions on $S \IH \tilde{M} \times \IB \times \IT$
\begin{eqnarray*}
\calf_{\IT} & = & \{ \; \{ (v, \beta , t) \; | \; t \leq t_0 \} \; | \; t_0 \in \IT \} \\
\calf_{\IT(\beta)} & = & \{ \; \{ (v, \beta , t) \; | \; t \leq t_0(\beta) \} \; | \; t_0\colon\IB \to \IT 
    \mbox{ a continuous function} \}  \\
\calf & = & \calf_{\IB} \cup \calf_{\IT( \beta )}.
\end{eqnarray*}
Observe that germs away from $\calf_{\IT}$ are the usual germs at infinity. Since the definition
of $\calf$ only involves the $\IB$- and the $\IT$-coordinate we will later use the same notation
for other (subsets of) spaces of the form $X \times \IB \times \IT$. 
We can reformulate the Foliated Control Theorem as follows.

\begin{proposition} \label{first-reduction}
The map in Theorem~\ref{fct} induces an  equivalence in $K$-theory if and only if the map
\[
\calc^{\Gamma} ( (S \IH \tilde{M} \times \IB \times \IT)_\angle , \cale_s )^{> \calf}
\to 
\calc^{\Gamma} ( (S \IH \tilde{M} \times \IB \times \IT)_\angle , \cale_w )^{> \calf}
\]
induces an equivalence in $K$-theory.
\end{proposition}
Observe that there is no longer a  compactness condition in the $\IB$-direction but instead 
of germs at infinity we now have germs away from $\calf$.
\begin{proof}
For the purpose of this proof we introduce the following abbreviation. For object support 
conditions $\calf^{\prime}$ and $\calf^{\prime \prime}$ and a morphism support 
condition $\cale$ on $S \IH \tilde{M} \times \IB \times \IT$ set
\[
(\cale, \calf^{\prime})^{> \calf^{\prime \prime}} = 
\calc^{\Gamma} ( (S \IH \tilde{M} \times \IB \times \IT)_\angle , \cale , 
                                   \calf^{\prime} )^{>\calf^{\prime \prime}} .
\]
We have the following commutative diagram
\[
\xymatrix{
(\cale_s , \calf_{\IB} )^{ > \calf_{\IT ( \beta )}}  \ar[r] \ar[d] &
(\cale_s )^{ > \calf_{\IT ( \beta )}}  \ar[r] \ar[d] &
(\cale_s )^{ > \calf}  \ar[d]  \\
(\cale_w , \calf_{\IB} )^{ > \calf_{\IT ( \beta )}}  \ar[r] & 
(\cale_w )^{ > \calf_{\IT ( \beta )}}   \ar[r]  &
(\cale_w )^{ > \calf}.   
          }
\]
Here the vertical map on the left is the map in the Foliated Control Theorem~\ref{fct} 
because under the presence of the $\calf_{\IB}$-object support condition there is no 
difference between germs away from $\calf_{\IT}$ (alias germs at infinity)
and germs away from $\calf_{\IT(\beta)}$.
According to Lemma~\ref{fibration-lemma}~\ref{fibration-lemma-drei} both rows yield fibration sequences in $K$-theory and 
it hence suffices to show that both categories in the middle admit an Eilenberg swindle and 
are therefore contractible.
Pick $\delta_t \in \cals$. Then a swindle is induced in both cases 
by the map $(x,\beta,t) \mapsto (x,\beta+\delta_t,t)$ on 
$S \IH \tilde{M} \times \IB \times \IT$,
compare \cite[Proposition~4.4]{Bartels-Farrell-Jones-Reich(topology)}. 
\end{proof}


\subsection{Second reduction -- Discretization}

We would like to use the cell structures on the chunks 
$|L_{i,\beta}| \times [\beta , \beta +1] \times [ t_i , t_{i+1} ]$. 
But we do not know how around the boundary of the squares 
$ [\beta , \beta +1] \times [ t_i , t_{i+1} ]$ the different cell structures
fit together. To avoid this problem we use a Mayer-Vietoris argument.

We define the following subsets of $\IB$ respectively $\IT$.
\begin{eqnarray*}
\IB_{e} = \bigcup_{ \stackrel{\beta \in \IN_0}{\beta \mbox{ \small even}}}
\left[ \beta , \beta+1 \right]  
& \quad &
\IB_{o} = \bigcup_{ \stackrel{\beta \in \IN_0}{\beta \mbox{ \small odd}}}
\left[ \beta , \beta+1 \right]   \\
\IT_{e} = \bigcup_{ \stackrel{i \in \IN}{i \mbox{ \small even}}}
\left[ t_i , t_{i +1} \right]  
& \quad &
\IT_{o} = \bigcup_{ \stackrel{i \in \IN}{i \mbox{ \small odd}}}
\left[ t_i , t_{i +1} \right]   .
\end{eqnarray*}
Moreover we set
\[
\IB_{e \cap o} = \IB_{e} \cap \IB_{o} = \IN \subset \IB 
\]
and 
\[
\IT_{e \cap o} = \IT_{e} \cap \IT_{o} = \{ t_1, t_2 , \dots \} \subset \IT   .
\]
For a subspace $Y \subset S \IH \tilde{ M } \times \IB \times \IT$
we denote the intersection of $Y$ with $(S \IH \tilde{ M } \times \IB \times \IT)_\angle$ by $Y_\angle$. 

The condition that there are no nontrivial morphisms between different path components of $Y$ 
will be denoted by $\cale_{\pi_0 (Y)}$ or briefly $\cale_{\pi_0}$. 
Formally this can be defined as the pullback of the morphism control condition consisting 
only of the diagonal on $\pi_0 ( Y )$ via the natural projection $Y \to \pi_0 ( Y )$. 
\begin{warning} \label{warning-components}
If $X$ is a subset of $Y$ one should not confuse $\cale_{\pi_0 (Y)}$ restricted to $X$, which 
is again denoted  $\cale_{\pi_0 ( Y )}$, with $\cale_{\pi_0 (X)}$.
\end{warning}

\begin{proposition} \label{second-reduction}
If for all 9 spaces $S \IH \tilde{M} \times \IB_p \times \IT_q$ with $p$, 
$q \in \{ e, o , e \cap o \}$
the maps
\[
\calc^{\Gamma} ( (S \IH \tilde{M} \times \IB_p \times \IT_q)_\angle , 
                        \cale_s \cap \cale_{\pi_0} )^{> \calf}
\to 
\calc^{\Gamma} ( (S \IH \tilde{M} \times \IB_p \times \IT_q)_\angle , 
                        \cale_w \cap \cale_{\pi_0} )^{> \calf}
\]
induce equivalences in $K$-theory then so does the map in Proposition~\ref{first-reduction}.
\end{proposition}

\begin{proof}
If one drops the extra $\cale_{\pi_0}$-condition on both sides,
this follows easily by applying a Mayer-Vietoris argument (compare Remark~\ref{Mayer-Viet}) 
in the $\IT$- and then again in the $\IB$-direction. 
However, dropping the $\cale_{\pi_0}$-condition does not
change the categories: if the support of a morphism in 
$\calc^{\Gamma} ( (S \IH \tilde{M} \times \IB_p \times \IT_q)_\angle ,\cale_s )$
or $\calc^{\Gamma} ( (S \IH \tilde{M} \times \IB_p \times \IT_q)_\angle ,\cale_w )$
violates this condition, it does so only on a set in $\calf$ and this can  be ignored 
since  we take germs away from $\calf$. For the $\IT$-direction this 
follows from 
the fact that in the definition of foliated and metric control we always require
a bound in the $\IT$-direction and the fact that the distance between the $t_i$ increases with $i$.
For the
$\IB$-direction it follows from the $\delta_t$-control in this direction.   
\end{proof}

From now on we will restrict our attention to the space 
$S \IH \tilde{M} \times \IB_e \times \IT_e$, all the other cases are completely analogous.


\subsection{Induction over the skeleta}

We next define a filtration for $Y=S \IH \tilde{M} \times \IB_e \times \IT_e$. Recall that
\[
S = \bigcup_{\beta , i} |L_{\beta , i}| \times [ \beta , \beta+1 ] \times [ t_i, t_{i+1} ] 
                                     \subset S \IH \tilde{M} \times \IB \times \IT
\]
and set 
\begin{eqnarray*} 
Y^{(-1)} & = &  (\overline{S \IH \tilde{M} \times \IB_e \times \IT_e )
  - S}  \\ 
Y^{(k)} & = & Y^{(-1)} \cup \bigcup_{\beta , i 
  \mbox{ even}} | L_{\beta , i}^{(k)} | \times [ \beta , \beta +1 ] \times [t_i, t_{i+1} ] .
\end{eqnarray*}

For a subspace $X \subset Y$ and a morphism support condition $\cale$ on $Y$ we define the 
object support condition
\[
X^\cale = \{ X^E | E \in \cale \} .
\]
(Recall that  $X^E=\{ y \in Y| \mbox{ there exists an }x \in X \mbox{ with } (x,y) \in E \}$
denotes the $E$-thickening of $X$ in $Y$.)

\begin{proposition}[Induction Step] \label{induction-step} 
For $k=0,1 , \dots ,n=\dim S \IH \tilde{ M } $ the relax control map 
\[
\xymatrix{
\calc^{\Gamma} ( Y^{(k)}_\angle , 
  \cale_s \cap \cale_{\pi_0 (Y) } )^{ > \calf \cup \left( Y^{(k-1)} \right)^{ \cale_s \cap \cale_{ \pi_0 (Y)}}}  
\ar[d] \\
\calc^{\Gamma} ( Y^{(k)}_\angle , 
  \cale_w \cap \cale_{\pi_0 (Y) } )^{ > \calf \cup \left( Y^{(k-1)} \right)^{ \cale_w \cap \cale_{ \pi_0 (Y)}}}  
         }
\]
induces an equivalence in $K$-theory.
\end{proposition}

This proposition is proven by
combining Proposition~\ref{foliated-excision}, Lemma~\ref{compare-euclidean} and Lemma~\ref{euclidean-standard}.
Before we proceed we note that Proposition~\ref{induction-step} implies the Foliated Control 
Theorem~\ref{fct}.

\begin{corollary} \label{corollary-induction-step} 
The map 
\[
\calc^{\Gamma} ( Y_\angle , \cale_s \cap \cale_{\pi_0 (Y)} )^{ > \calf}
\to 
\calc^{\Gamma} ( Y_\angle , \cale_w \cap \cale_{\pi_0 (Y)} )^{ > \calf}
\]
induces an equivalence in $K$-theory.
By Propositions~\ref{second-reduction} and \ref{first-reduction} this implies the Foliated 
Control Theorem~\ref{fct}.
\end{corollary}

\begin{proof}[Proof of the Corollary]
Let us abbreviate
\begin{eqnarray*} 
\calc^{(k)}_x  & = & \calc^{\Gamma} ( Y^{(k)}_\angle , \cale_x \cap \cale_{\pi_0 (Y)} )^{ > \calf} \\
\calc^{(k,k-1)}_x & = & \calc^{\Gamma} ( Y^{(k)}_\angle , \cale_x \cap \cale_{\pi_0 (Y) })^{ > \calf \cup \left( Y^{(k-1)} \right)
                                       ^{ \cale_x \cap \cale_{ \pi_0 (Y)}}} 
\end{eqnarray*}
for $x = w,s$. By definition of $\cale_s$ and $\cale_w$ we have 
$\calc^{(-1)}_w = \calc^{(-1)}_s$. 
By Lemma~\ref{fibration-lemma}~\ref{fibration-lemma-drei} and  
Remark~\ref{proper-remark}
the sequence 
\[
\xymatrix
{
\calc^{(k-1)}_x \ar[r]  &
\calc^{(k)}_x \ar[r]  &
\calc^{(k,k-1)}_x
}
\]
induces for $x=w$ or $x=s$ a fibration sequence in $K$-theory and relaxing 
control from $s$ to $w$ yields a map of fibration sequences.
The result follows by induction using \ref{induction-step}.
\end{proof}

Observe that in Proposition~\ref{induction-step} taking germs away from $Y^{(k-1)}$
means in particular that everything that is relevant happens over $S$, i.e.\ the region covered by 
the cell structures (compare (\ref{definition-S})). Hence we can 
ignore the difference between strong and very strong control. (Formally this is an application of 
Lemma~\ref{AFfequivalence} from the Appendix.) 
Moreover we can drop
the ${\angle}$-subscript because 
\[
Y^{(k)} - Y^{(k)}_\angle \subset Y^{(-1)} \subset Y^{(k-1)}.
\]
We obtain the following lemma.
\begin{lemma}
The map in Proposition~\ref{induction-step} is equal  to the map
\[
\xymatrix{
\calc^{\Gamma} ( Y^{(k)} , \cale_{vs} \cap \cale_{\pi_0 ( Y)})
      ^{> \calf \cup ( Y^{(k-1)} )^{\cale_{vs} \cap \cale_{\pi_0(Y)}}} \ar[d] \\
\calc^{\Gamma} ( Y^{(k)} , \cale_{w}  \cap \cale_{\pi_0 ( Y)})
      ^{> \calf \cup ( Y^{(k-1)} )^{\cale_{w} \cap \cale_{\pi_0(Y)}}}.
         }
\]
\end{lemma} 

Recall that $L^{ \{ k \} }_{ \beta ,i }$ denotes 
those $k$-cells in 
$L_{\beta , i }$ which do not lie in the $(-1)$-st filtration step, i.e.\ in $Y^{(-1)}$.
Set 
\begin{eqnarray*}
Z^{(k)} 
& = & 
\coprod_{\stackrel{e \in L_{\beta , i}^{ \{ k \} }}{\beta ,i \mbox{ {\tiny even}} }} (A_e \times B_e)  
              \times [ \beta , \beta+1 ] \times [ t_i , t_{i+1} ] \\
\partial Z^{(k)} 
& = & 
\coprod_{\stackrel{e \in L_{\beta , i}^{ \{ k \} }}{\beta ,i \mbox{ {\tiny even}} }}  
                           \partial(A_e \times B_e)  \times [ \beta , \beta+1 ] \times [ t_i , t_{i+1} ] .
\end{eqnarray*}
There is a natural map 
\[
g \colon Z^{(k)} \to Y^{(k)}
\]
induced by the charts $g_e$. The map induces a homeomorphism 
$(Z^{(k)} - \partial Z^{(k)}) \to (Y^{(k)} - Y^{(k-1)})$. 
Recall that the object
support condition $\calf$ is defined for every space with $\IB$- and $\IT$-coordinate, 
compare the beginning of Subsection~\ref{subsection-first-reduction}.
In particular $g^{-1} \calf$ will again be denoted $\calf$.
The following proposition is the crucial step in the proof of the Foliated Control Theorem and 
should be thought of as a ``foliated excision'' result. 
It allows to separate cells. (Note the $\cale_{\pi_0(Z^{(k)})}$-condition in the source.)

\begin{proposition}[Excision of the $(k-1)$-skeleton] \label{foliated-excision}
Let $\cale$ denote either $\cale_{vs}$ or $\cale_w$. In both cases the natural map 
$g \colon Z^{(k)} \to Y^{(k)}$ induces an equivalence of categories
\[
\xymatrix{
\calc^{\Gamma} (Z^{(k)} , g^{-1} \cale \cap \cale_{\pi_0(Z^{(k)} )} )
^{ > \calf \cup \left( \partial Z^{(k)} \right)^{ g^{-1} \cale \cap \cale_{\pi_0 ( Z^{(k)} ) } } } 
\ar[d] \\
\calc^{\Gamma} (Y^{(k)} , \cale \cap \cale_{\pi_0( Y )} )
^{> \calf \cup \left( Y^{(k-1)}  \right)^{ \cale \cap \cale_{\pi_0 ( Y ) } } }
         }
\]
for $k=0 , 1, \dots , n$.
\end{proposition}

\begin{proof}
Observe that the map factorizes over 
\[
\calc^{\Gamma} (Z^{(k)} , g^{-1} \cale \cap g^{-1} \cale_{\pi_0( Y ) })^{ >\calf \cup 
\left( \partial Z^{(k)} \right)^{ g^{-1} \cale \cap g^{-1}\cale_{\pi_0 ( Y ) } } } .
\]
Since $g \colon Z^{(k)} \to Y^{(k)}$ induces a homeomorphism 
$Z^{(k)} - \partial Z^{(k)} \to Y^{(k)} -  Y^{(k-1)}$ 
and all conditions are simply pulled back along $g$
the second map in this factorization clearly induces an 
equivalence. 
We see that the crucial point is whether the forget 
control map from $g^{-1} \cale \cap \cale_{\pi_0( Z^{(k)} )}$-control to 
$g^{-1} \cale \cap g^{-1} \cale_{\pi_0( Y ) }$-control induces an equivalence.
Note that $\cale_{\pi_0 ( Z^{(k)} )}$ does not allow morphisms between different cells and 
is hence a lot stronger than $g^{-1} \cale_{\pi_0( Y)}$ which only separates the different $[ \beta , \beta +1] \times [t_i , t_{i+1} ]$-blocks.
Formally the result will 
be a consequence of Lemma~\ref{AFfequivalence} in the Appendix. 
We will apply that lemma to the case 
$X=Z^{(k)}$, $A= \partial Z^{(k)}$,
$\cale^{\prime }= g^{-1}\cale \cap \cale_{\pi_0( Z^{(k)} )}$, 
$\cale^{\prime \prime}= g^{-1}\cale \cap  g^{-1}\cale_{\pi_0( Y )}$, 
$\calf_0= \emptyset$ and $\calf=\calf$. 
We only formulate the argument for $\cale=\cale_w$. The $\cale_{vs}$-case is easier and can be 
obtained by setting $\alpha=0$, compare Remark~\ref{alpha-zero-control}.

Let $\alpha \in \IN_0$ and $\delta_t \in \cals$ be given.
Write $E_{\alpha , \delta_t }$ for the subset $E \in \cale$ determined
as in Definition~\ref{fol-control-w-decay}~(ii) 
by $\alpha$ and $\delta_t$ (we suppress the $t_0$). 
As in Proposition~\ref{S-properties}~(ii)
define for the given $\delta_t=\step_{\Jt} ( (\delta_i)) \in \cals$ the function 
$\Delta_t = \step_{\Jt} ( (\Delta_i(\delta_i)))$ for 
$t \geq t_{i_0}$. In particular, $\delta_i < \e_i$ for all $i \geq i_0$.
By Proposition~\ref{S-properties}~(ii) we have 
$E_{\alpha, \Delta_t} \in \cale$. 
Using our usual coordinates $\beta \in \IB$, $t \in \IT$ 
we define $F$ as the union of
\begin{eqnarray*}
F_{t \leq t_{i_0}} & = & \{ (v,\beta,t) \; | \; t \leq t_{i_0} \}, \\
F_{\beta \geq i } & = & S \IH \tilde{ M} \times 
       \bigcup_{ \beta \geq i } [ \beta , \beta +1 ] \times [ t_i , t_{i+1} ] \qquad \mbox{and}  \\
F_{\beta \leq \alpha} & = & \{ (v,\beta,t) \; | \; \beta \leq \alpha \}.
\end{eqnarray*}
Observe that $F_{t \leq t_{i_0}} \cup F_{\beta \geq i} \in \calf_{\IT ( \beta )}$, 
$F_{\beta \leq \alpha } \in \calf_{\IB}$ and hence $F \in \calf$.
Let $E'' \subset g^{-1}E_{\alpha,\delta_t}$ resp. $E' \subset g^{-1}E_{\alpha,\Delta_t}$ 
denote the subset
of all pairs of points that satisfy in 
addition the $g^{-1}\cale_{\pi_0(Y)}$- resp. $\cale_{\pi_0(Z^{(k)})}$-condition. We need to
check that with this notation the condition in Lemma~\ref{AFfequivalence} is satisfied. 
Since $E_{\alpha,\delta_t} \subset E_{\alpha,\Delta_t}$ it suffices to show that if
$ (x,x') \in E''$ and 
$x \notin (\partial Z^{(k)})^{E'} \cup F$ then
$(x,x')$ satisfies the $\cale_{\pi_0(Z^{(k)})}$-condition:
let $(y,\gamma,t) = g(x)$ and $(y',\gamma',t') = g(x')$ (here $\gamma$ is the $\IB$- and $t$ the $\IT$-coordinate).
There are $i,\beta$ such that
$t_i \leq t,t' \leq t_{i+1}$ and $\beta \leq \gamma,\gamma' \leq \beta+1$, 
since $(x,x')$ satisfy the $\cale_{\pi_0(Y)}$-condition. 
We know $i \leq \beta$   (since $x \notin F_{\beta \geq i}$), 
        $\beta > \alpha$ (since $x \notin F_{\beta \leq \alpha}$)
    and $i > i_0$ (since $x \notin F_{t \leq t_{i_0}}$). 
In particular, $\delta_i < \e_i \leq \e_{L_{\beta,i}}$. 
There is $e \in L^{ \{k\} }_{\beta,i}$ such that  
$y \in e - \partial e^{(\alpha,\Delta)}$     
where $\Delta = \Delta_i(\delta_i) \geq \Delta_{L_{\beta,i}}(\alpha,\delta_i)$, 
since $x \notin (\partial Z^{(k)})^{E'}$. 
Finally, $d_{F_{geo}}(y,y') \leq (\alpha,\delta_i)$,
since $(x,x')$ satisfies the $E_{\alpha,\delta_t}$-condition. 
All this allows the application of Lemma~\ref{flow-cell-structure-constants}~\ref{fcsc-eins} 
(compare also Lemma~\ref{simultaneous-for-all}~\ref{simultaneous-for-all-Delta})
to conclude that $y' \in e$. 
Thus $(x,x')$ does indeed satisfy the $\cale_{\pi_0(Z^{(k)})}$-condition. 
\end{proof}


\subsection{Comparison to a Euclidean standard situation}

According to the last proposition we can assume that all the cells which are new in the $k$-th 
step do not talk to each other through 
morphisms (this is formalized in the $\cale_{\pi_0 (Z^{(k)})}$-condition).
In the next step we will use the fact that each flow cell also has a security zone around it 
on which we have a very precise control over the foliation and the metric to prove that
the situation is equivalent to a Euclidean standard situation.
More precisely define 
\begin{eqnarray*}
W^{(k)} & = & 
  \coprod_{\stackrel{ e \in L_{\beta , i}^{ \{ k \} }}{\beta,i \mbox{ {\tiny  even}}} } 
     \left( \IR \times \IR^{n-1} \right) \times [ \beta , \beta +1] \times [ t_i , t_{i+1} ],  
\end{eqnarray*}
We equip each $(\IR \times \IR^{n-1}) \times [\beta , \beta+1]$ 
with the standard Euclidean metric and with the 
foliation $F_{\IR^n}$ (i.e.\ the foliation by lines parallel to the $\IR$-factor). 
We define in the obvious way 
the foliated and the metric control structure on $W^{(k)}$.
Namely in both cases we impose $\cale_{\pi_0 ( W^{(k)} )}$-control, 
i.e.\ different components are infinitely far apart.
We then let $\cale_{met}$ denote metric control with decay speed $\cals$ together with 
$\cale_{\pi_0( W^{(k)} )}$ and $\cale_{fol}$ the foliated control with decay speed $\cals$ 
together with $\cale_{\pi_0 ( W^{(k)} )}$.
We consider $Z^{(k)}$ as a subset of $W^{(k)}$ and denote this inclusion by $i \colon Z^{(k)} \to W^{(k)}$. 
Note that this map induces an equivalence on 
$\pi_0$ and hence $i^{-1} \cale_{\pi_0(W^{(k)})}= \cale_{\pi_0 ( Z^{(k)} )}$. Now consider 
the following situation. 
\[
\xymatrix{ 
( Y^{(k)} , \cale_{vs}) \ar[d] & Z^{ (k) } \ar[d]^{=} \ar[l]_-g \ar[r]^-i  & ( W^{(k)} , \cale_{met} ) \ar[d]  \\
( Y^{(k)} , \cale_{w}) & Z^{ (k) } \ar[l]_-g \ar[r]^-i  & ( W^{(k)} , \cale_{fol} )  .
}
\]
The following proposition reduces our problem to the  Euclidean standard situation on the right of the diagram above.

\begin{lemma}[Comparison to a Euclidean situation]  \label{compare-euclidean}
For $k=0, 1, \dots ,n$ the categories
\[
\calc^{\Gamma} ( Z^{ (k)} , \cale )^{> \calf \cup ( \partial Z^{(k)} )^{ \cale }}
\]
with $\cale= g^{-1} \cale_{vs} \cap \cale_{\pi_0 ( Z^{(k)} )}$ resp.\ 
$\cale=i^{-1} \cale_{met}$ are equivalent.
The same holds for the pair of control conditions 
$\cale= g^{-1} \cale_w \cap \cale_{\pi_0 ( Z^{(k)} )}$ and 
$\cale=i^{-1} \cale_{fol}$.
\end{lemma}

\begin{proof}
We need to 
show that the four horizontal relax control maps in the 
following diagram are equivalences.
\[
\xymatrix{
g^{-1} \cale_{vs} \cap \cale_{\pi_0 ( Z^{(k)} ) } \ar[d] & 
g^{-1} \cale_{vs} \cap \cale_{\pi_0 ( Z^{(k)} ) } \cap i^{-1} \cale_{met} \ar[l] \ar[r] \ar[d] &
i^{-1} \cale_{met} \ar[d] \\
g^{-1} \cale_{w} \cap \cale_{\pi_0 ( Z^{(k)} ) } & 
g^{-1} \cale_{w} \cap \cale_{\pi_0 ( Z^{(k)} ) } \cap i^{-1} \cale_{fol} \ar[l] \ar[r] &
i^{-1} \cale_{fol}. 
         }
\]
We only treat the lower left hand horizontal map. The other cases are analogous.
Formally the argument is an application of Lemma~\ref{AFfequivalence}.

Let $\alpha$ and $\delta_t \in \cals$ be given.
Let $E''$ be determined by $g^{-1}E_{\alpha, \delta_t}$ and the 
$\cale_{\pi_0 (Z^{(k)})}$-condition.
Define $E'$ by $i^{-1} E_{\alpha , C_t  \delta_t } \cap E''$, where $C_t$ stems from
Proposition~\ref{S-properties}~\ref{S-properties-drei} which also tells us that $C_t \delta_t \in \cals$.
Choose $i_0$ such that for all $i \geq i_0$ we have 
$C_i \delta_i \leq \epsilon_i$, see Lemma~\ref{construct-T} \ref{cT-vier} and \ref{cT-drei}.
Define $F_{t \leq t_{i_0} }$, $F_{\beta \geq i}$ and $F_{\beta \leq (\mu_n /5) \cdot  \alpha}$ as in the
proof of Proposition~\ref{foliated-excision} (but note the constant $\mu_n /5$) and let $F$ be the union
of these three sets. The condition in Lemma~\ref{AFfequivalence} is implied if one can show that
two points in the same cell of a cell structure $L_{\beta,i}$ which are less than 
$(\alpha, \delta_{t_i})$-foliated apart, when measured inside the manifold 
are $(\alpha , C_{t_i} \delta_{t_i})$-controlled when measured in Euclidean space using the charts. 
At least this should be true
away from the set $F$. But this is the content of Lemma~\ref{simultaneous-for-all}~\ref{simultaneous-for-all-C} which says 
that
Remark~\ref{compare-foliated-distance} applies with $C_i$ instead of $C_e$ and with $\epsilon_i$ 
instead of $\epsilon_e$ if $\beta \leq i$.
Note that we can assume 
$\alpha \leq 5 \mu_n^{-1} \beta$ and $t \geq t_{i_0}$ which translates into $\delta_i \leq C_i^{-1} \epsilon_i$.
These are just the assumptions in Remark~\ref{compare-foliated-distance}.
\end{proof}

It remains to prove the comparison result for the Euclidean standard situation.

\begin{lemma}  \label{euclidean-standard}
For $k=0,1,\dots ,n$ the map 
\[
\xymatrix{
\calc^{\Gamma} ( Z^{(k)} , i^{-1} \cale_{met}  )^{ > \calf 
\cup ( \partial Z^{(k)} )^{i^{-1} \cale_{met} } } \ar[d] \\
\calc^{\Gamma} ( Z^{(k)} , i^{-1} \cale_{fol}  )^{ > \calf 
\cup ( \partial Z^{(k)} )^{i^{-1} \cale_{fol} } } 
} 
\]
induces an equivalence in $K$-theory.
\end{lemma}

\begin{proof}
Write $Z^{(k)} = Z_t^{(k)} \cup Z_l^{(k)}$ where the first subspace uses only transversal cells
and the second only long cells. It suffices to check the claim separately for $Z_t^{(k)}$ and  
$Z_l^{(k)}$. On transversal cells $i^{-1}\cale_{met}$ and $i^{-1}\cale_{fol}$ 
and the relevant thickenings of the boundaries agree 
and we are done for 
$Z_t^{(k)}$. If $A_e = [a,b]$ let $\dd_-^e = \{ a \} \x B_e$ and write 
$\dd(A_e \x B_e) = \dd_-^e \cup \dd_+^e$ where $\dd_+^e \cap \dd_-^e = \{ a \} \x \dd B_e$.
Set
\begin{eqnarray*}
\dd_\pm Z_t^{(k)} &  = &   
   \coprod_{\stackrel{ e \mbox{ \tiny long cell in } L_{\beta , i}^{ \{ k \} }   }{\beta,i \mbox{ {\tiny  even}}} } 
      \dd^e_\pm  \x [ \beta , \beta +1] \times [ t_i , t_{i+1} ].
\end{eqnarray*}   
Let $\cale$ denote either $i^{-1} \cale_{met}$ or $i^{-1} \cale_{fol}$, then according to 
Lemma~\ref{fibration-lemma}~\ref{fibration-lemma-drei} the $K$-theory of the categories we are interested in is the cofiber
of the map induced by
\[
\xymatrix{
\calc^{\Gamma}( Z_l^{(k)} , \cale , ( \partial_- Z_l^{(k)} )^{\cale} )^{ > \calf \cup (\partial_+ Z_l^{(k)})^{\cale} }
\ar[r] &
\calc^{\Gamma}( Z_l^{(k)} , \cale )^{ > \calf \cup (\partial_+ Z_l^{(k)})^{\cale} }.
         }
\]
For both choices of $\cale$ the map 
$(v,\gamma,t) \mapsto (v + (\delta_t,0),\gamma,t)$ with some fixed $\delta_t \in \cals$
induces an  Eilenberg swindle for the category on the right. (Here $v=(a,b)$ refers to the coordinates in  $A_e \times B_e$.)
By Remark~\ref{proper-remark} the category on the left can be identified with 
\[
\calc^{\Gamma}( \partial_-  Z_l^{(k)} , \cale )^{ > \calf \cup (\partial_+ Z_l^{(k)})^{\cale} }.
\]
Note that $\partial_- Z_l^{(k)}$ consists only of transversal pieces
and we can repeat the argument from the beginning.
\end{proof}

This finishes the proof of Proposition~\ref{induction-step} and hence by 
Corollary~\ref{corollary-induction-step} the proof of the Foliated Control Theorem~\ref{fct}.


\typeout{---------------------collapsing-------------------------}

\section{From strong control to continuous control}
\label{section-collapsing}

In this section we explain the map (8) in the main diagram in Section~\ref{section-proof} which connects up the ``strong'' control
($\cale_s$-control)
we obtained so far with the equivariant continuous control on the space $X( \infty ) \times \IT$ 
(denoted $\cale_{\Gamma cc} ( X(\infty) )$-control).


\subsection{A space with infinite cyclic isotropy}
\label{subsection-construct-X}
We recall the construction of the metric space $X$ and the map $p_X : S \IH \tilde{M} \to X$ from 
\cite[Section~14]{Bartels-Farrell-Jones-Reich(topology)}. 
One can collapse $\IH \tilde{M}$ to the $\Gamma$-compact space $\IH_{ [-1 , 1] } \tilde{M}$
by projecting $\IH_{(-\infty , -1]} \tilde{M}$ in the obvious way to $\IH_{\{ -1 \}} \tilde{M}$ and  likewise
$\IH_{[ 1, \infty )} \tilde{M}$ to $\IH_{\{ 1 \} } \tilde{M}$. Similarly $X$ is obtained from
$S \IH \tilde{M}$, where additionally the fibers of the bundle $S \IH \tilde{M} \to \IH \tilde{M}$
over $\IH_{\{-1\}} \tilde{M}$ and $\IH_{\{ 1 \} } \tilde{M}$ are collapsed to points. This collapsing
map is the map $p_X : S \IH \tilde{M} \to X$. 
More details of the construction can be found  before Proposition~14.5 in \cite{Bartels-Farrell-Jones-Reich(topology)}.
We can also project all the way down to $\tilde{M}$ and we hence obtain
a factorization of the natural projection $S \IH \tilde{M} \to \tilde{M}$ over $X$.
Restricted to $S \IH_{ [-0.5 , 0.5] } \tilde{M}$ the map
$p_X$ is essentially the identity and we can hence consider  $S \IH_{ [-0.5 , 0.5] } \tilde{M}$ and for every $\beta \geq 0$ also 
$S \IH \tilde{M}^{\leq \mu_n \beta}$ as a subset of $X$. (Recall that $S\IH \tilde{M}^{\leq \mu_n \beta}$ denotes the union of all
leaves which are shorter than $\mu_n \beta$.)
Later the metric 
properties of $p_X$ will be important: $p_X$ does not increase distances and for $i$ large the map $p_X$ contracts 
$S \IH \tilde{M} -  S \IH_{[-i,i]} \tilde{M}$ rather strongly. The precise statement is  \cite[14.5]{Bartels-Farrell-Jones-Reich(topology)}.

We define $X(\beta)$ via the following pushout diagram whose horizontal arrows are $\Gamma$-cofibrations
\[
\xymatrix
{
S \IH \tilde{M}^{\leq \mu_n \beta} \ar[d] \ar[r]  &
X \ar[d]^-{p_{\beta}} 
\\
\pi_0( S \IH \tilde{M}^{\leq \mu_n \beta}   )  \ar[r] &
X(\beta) .
}
\] 
Since $S \IH \tilde{M}^{\leq \mu_n \beta} \subset S \IH \tilde{M}^{\leq \mu_n (\beta +1) }$ we obtain natural maps
$c_{\beta} \colon X( \beta ) \to X ( \beta +1 )$ and we define $X( \infty )$ as the 
mapping telescope model for $\hocolim_{ \beta \geq 0} X(\beta)$, i.e.\
as the coequalizer of 
\[
\xymatrix{
\coprod_{\beta \geq 0} X(\beta) \ar@<0.5ex>[r] \ar@<-0.5ex>[r] & \coprod_{\beta \geq 0} X(\beta) \x [\beta,\beta+1],\\
         }
\]  
where the maps are given by sending $x \in X( \beta )$ to $(x,\beta +1)$ respectively to $(c_{\beta } (x), \beta )$.
Note that $X(\infty)$ is a $\Gamma$-space all whose isotropy groups
are trivial or infinite cyclic. 
We obtain natural maps
\[
S\IH \tilde{M} \times \IB  \xrightarrow{q} X \times \IB  \xrightarrow{p}  X(\infty),
\]
where $q = p_X \times \id_{\IB}$ and $p$ is induced from the maps $p_{\beta}$.


\subsection{Strong control maps to continuous control}
 
We will now check that $q$ does in fact define the arrow labeled (8) in the main diagram in Section~\ref{section-proof}.

\begin{proposition} \label{pro-strong-maps-to-continous}
The map $q$ induces a functor 
\[
\calc^{\Gamma}( (S \IH \tilde{M} \times \IB \times \IT)_\angle , \cale_s, \calf_{\IB} )^{\infty}
\to
\calc^{\Gamma} ( X \times \IB \times \IT , (p \x \id_\IT)^{-1} \cale_{\Gamma cc}(X(\infty)), \calf_{\IB} )^{\infty} .
\]
\end{proposition}

\begin{proof}
Note that the restriction of $q\times \id_{\IT}$ to $(S \IH \tilde{M} \times \IB \times \IT)_\angle$ is a proper map.
Compare the discussion of the map (4) in the outline of the proof in Section~\ref{section-proof}.
We need to check that  $q$ maps $\cale_s$-control to $p^{-1} \cale_{\Gamma cc}$-control.
This has roughly the following reasons: 
Firstly, since $q$ does not increase distances (\cite[14.5(i)]{Bartels-Farrell-Jones-Reich(topology)}) it maps very strong control ($\cale_{vs}$-control, compare
Definition~\ref{definition-weak-strong-control})
to continuous control. Since $\cale_{s}$-control  implies $\cale_{vs}$ over the set $S$ (which was defined at the beginning of 
Subsection~\ref{statement-fct}), $\cale_{s}$-control implies continuous control over $S$. 
Secondly, we need to deal with points on a geodesic $g$ that is collapsed to a point  
in $X(\infty)$. Here foliated control ($\cale_w$-control) already implies continuous control. 
Thirdly, we are left with points that are not in $S$ because they have a large $\IH$-coordinate. Here $\cale_s$-control
implies only bounded control, but $q$ contracts this part very strongly (\cite[14.5(ii)]{Bartels-Farrell-Jones-Reich(topology)}) and produces continuous control.
For a careful argument one needs to construct suitable 
neighborhoods of $g$ that are invariant under the stabilizer of $g$ (compare~\cite[15.2]{Bartels-Farrell-Jones-Reich(topology)}) and to use the fact that we did choose
the tubular neighborhoods $T_{\beta ,i}$ monotone decreasing and such that 
$\bigcap_{i \in \IN} T_{\beta,i} = S \IH M^{\leq \mu_n \beta}$.
\end{proof}


\typeout{---------------------appendix -------------------------}

\section{Appendix} \label{appendix}

In this Appendix we collect a couple of facts which are more easily treated independently from the context 
in which they were used in the main text. 

\subsection{Homotopy finite chain complexes} \label{appendix-homotopy-finite}

Let $\overline{\cala}$ be an additive category and $\cala \subset \overline{\cala}$ a full additive 
subcategory. We think of $\cala$ as the category of ``finite'' objects in $\overline{\cala}$, compare 
Subsection~\ref{subsection-homotopy-finite-chain-complexes}.
Given such a situation we denote by $\ch \overline{\cala}$ the category of bounded below chain complexes in $\overline{\cala}$.
The notion of chain homotopy leads to a notion of weak equivalence and we define cofibrations to be those chain maps
which are degreewise the inclusion of a direct summand. The category
$\ch \overline{\cala}$ becomes a Waldhausen category (a category with cofibrations and weak equivalences in the sense of \cite{Waldhausen(1126)}).
We define $\ch_f \cala$
as the the full subcategory of $\ch \overline{\cala}$ whose objects are bounded below and above complexes
where the object in each degree of the chain complex lies in $\cala$. 
Furthermore we define 
$\ch_{\hf} \cala$
to be the full subcategory of chain complexes in $\ch \overline{\cala}$ which are homotopy equivalent to a complex in
$\ch_{f} \cala$. The categories $\cala$, $\ch_f \cala$ and $\ch_{hf} \cala$ inherit a Waldhausen structure 
from $\ch \overline{ \cala }$. 
\begin{lemma} \label{lemma-inclusions-equivalences}
The natural inclusions
\[
\cala \to \ch_f \cala \to \ch_{\hf} \cala
\]
induce equivalences in connective $K$-theory.
\end{lemma}
\begin{proof}
For the first map see \cite{Brinkmann(1979)} and \cite{Thomason-Trobaugh(schemes)} or \cite{Cardenas-Pedersen(karoubi)}.
The second inclusion induces an equivalence by a standard application of the 
Approximation Theorem~1.6.7 in \cite{Waldhausen(1126)}.
(Mimic the mapping cylinder argument on page~380 in \cite{Waldhausen(1126)} for chain complexes.)
\end{proof}

\subsection{The tilde-construction} \label{appendix-tilde-construction}

In Section~\ref{section-transfer} we were forced to artificially enlarge the 
Waldhausen category $\ch_{\hf} \calc^{\Gamma} ( \tilde{E} \times \IT ; \cale )$
to a category with the same $K$-theory in order to define a transfer. In fact this enlargement
is most easily treated in the generality of Waldhausen categories. In this subsection we briefly describe how such an 
enlargement is formally defined and we explain that the natural inclusion defines an equivalence in $K$-theory
under mild conditions.

Let $\calw$ be a Waldhausen category.  We additionally assume:
\begin{enumerate}
\item[(M)]
Cofibrations in $\calw$ are monomorphisms.
\item[(H)]
There is a functor $\calc \to \cald$ to some other category such that precisely the 
weak equivalences are mapped to isomorphism in $\cald$.
\item[(Z)]
The category admits a cylinder functor and satisfies the cylinder axioms, compare page 348 in \cite{Waldhausen(1126)}.
\end{enumerate}
Note that (H) implies the 
saturation axiom, i.e.\ the ``two out of three''-axiom for weak equivalences (see \cite{Waldhausen(1126)} page 327).

Analogously to the Waldhausen categories $F_m \calw$ defined on page 324 in \cite{Waldhausen(1126)}
we define the Waldhausen category $F_{\infty} \calw$ whose objects 
\[
\xymatrix{
C= ( C_0 \ar@{^{(}->}[r]^-{j_0} & C_1 \ar@{^{(}->}[r]^-{j_1}  & C_2 \ar@{^{(}->}[r]^-{j_2}  & \dots )
         }
\]
are infinite sequences of cofibrations in $\calw$. We denote the full Waldhausen subcategory on those objects where all
the $j_i$ are additionally weak equivalences by $\widehat{\calw}$. The shift functor $\sh: \widehat{\calw} \to \widehat{\calw}$
sends $C_0 \to C_1 \to \dots$ to $C_1 \to C_2 \to \dots $, i.e.\ it simply forgets $C_0$. 
There is an obvious natural transformation 
$\tau$ from the identity functor to the shift functor.

We define $\widetilde{ \calw }$ to be the category whose objects are the same as the objects in $\widehat{ \calw}$
and where the set of morphisms between $C$ and $D$ is given by the colimit over the following commutative
diagram (which is indexed over the lattice points in a $\frac{3}{8}$-th plane).
\[
\xymatrix{
\dots & \dots &  \dots & \\
\mor_{\widehat{\calw}}( \sh C , D ) \ar[u]^-{\sh} \ar[r]^-{\tau_{\ast}} &
\mor_{\widehat{\calw}}( \sh C , \sh D ) \ar[u]^-{\sh} \ar[r]^-{\tau_{\ast}} &
\mor_{\widehat{\calw}}( \sh C , \sh^2 D ) \ar[u]^-{\sh} \ar[r]^-{\tau_{\ast}} & \dots \\
& \mor_{\widehat{\calw}}( C , D ) \ar[u]^-{\sh} \ar[r]^-{\tau_{\ast}} &
\mor_{\widehat{\calw}}( C , \sh D ) \ar[u]^-{\sh} \ar[r]^-{\tau_{\ast}} & \dots
         }
\]
Assumption (M) implies that all maps in this diagram are inclusions.
We define the cofibrations in $\widetilde{\calw}$ to be those morphisms which up to an isomorphism can be represented by
a cofibration in $\widehat{\calw}$ in the above colimit. Similarily $C \to D$ is a weak equivalence
if it can be represented by a weak equivalence $\sh^n C \to \sh^m D$ in $\widehat{\calw}$ for some $n$ and $m$.
A lengthy but straightforward argument shows that 
these structures indeed define the structure of a Waldhausen category on $\widetilde{\calw}$.
Assumption (H) is used to verify that isomorphisms are weak equivalences.

\begin{proposition} \label{enlarge-equivalence}
Suppose $\calw$ satisfies (M), (H) and (Z) then 
the natural inclusion 
\[
\calw \to \widetilde{\calw}
\]
which sends $C$ to $C \xrightarrow{=} C \xrightarrow{=} C \xrightarrow{=} \dots $
induces an equivalence on connective $K$-theory. 
\end{proposition} 
\begin{proof}
This is an application of Waldhausen's Approximation Theorem~1.6.7 in \cite{Waldhausen(1126)}.
\end{proof}

\subsection{Swan group actions on $K$-theory} \label{appendix-swan}

In this subsection we will briefly describe several versions of the Swan group and how it acts on $K$-theory.
We will use the notation
\[
\Sw ( \Gamma ; \IZ ) \quad \mbox{ and } \quad \Sw^{\fr} ( \Gamma ; \IZ ) 
\]
for the $K_0$-group of $\IZ\Gamma$-modules which are finitely generated as $\IZ$-modules respectively finitely generated 
free as $\IZ$-modules.
In both cases the relations are the additivity relation given for all (not necessarily $\IZ$- or $\IZ \Gamma$-split) exact sequences.
Furthermore we will need the ``chain complex version''
\[
\Sw^{\ch } ( \Gamma ; \IZ ),
\]
which is defined as the free abelian group on isomorphism classes of all bounded below complexes $C_{\bullet}$ of $\IZ \Gamma$-modules
satisfying
\begin{enumerate}
\item
The homology $H_{\ast} ( C_{\bullet} )$ is finitely generated 
as an abelian group (and in particular concentrated in finitely many degrees).
\item
The modules in each degree are free as $\IZ$-modules.
\end{enumerate}
modulo the relations
\begin{enumerate}
\item
A short exact sequence $0 \to C_{\bullet} \to D_{\bullet} \to E_{\bullet} \to 0$ of $\IZ \Gamma$-chain complexes yields
$[ D_{\bullet} ] = [ C_{\bullet} ] + [ E_{\bullet} ]$.
\item
If $C_{\bullet} \to D_{\bullet}$ is a $\IZ \Gamma$-chain map which induces an isomorphism on homology then $[C_{\bullet}] = [D_{\bullet}]$.
\end{enumerate} 

There is a natural map  $i:\Sw^{\fr} ( \Gamma ; \IZ ) \to \Sw ( \Gamma ; \IZ )$ and considering a module
as a chain complex concentrated in degree $0$ defines a map
$j: \Sw^{fr} ( \Gamma ; \IZ ) \to \Sw^{\ch} ( \Gamma ; \IZ )$. More remarkable is the map
\begin{eqnarray*}
\chi : \Sw^{\ch} ( \Gamma ; \IZ ) & \to & \Sw ( \Gamma ; \IZ ) \\
C_{\bullet} & \mapsto & \sum (-1)^i [ H_i ( C_{\bullet} ) ] .
\end{eqnarray*}
\begin{proposition} \label{all-three-isos}
All three maps in the commutative diagram
\[
\xymatrix{
\Sw^{\fr} ( \Gamma ; \IZ ) \ar[r]^-j \ar@/_3ex/[rr]_-i &
\Sw^{\ch} ( \Gamma ; \IZ ) \ar[r]^-{\chi} &
\Sw ( \Gamma ; \IZ ) 
         }
\]
are isomorphisms.
\end{proposition}
\begin{proof} For a $\IZ \Gamma$-module $M$ we denote by $TM$ the $\IZ$-torsion submodule. For $[M] \in \Sw (\Gamma ; \IZ)$
there exists a $\IZ \Gamma$-resolution $F_{\bullet}(TM)=(F_1 \to F_0)$ of $TM$ by modules which are finitely generated free as 
$\IZ$-modules. The map $\Phi : [M] \mapsto [F_1] - [F_0] + [M/TM]$ is a well defined inverse for $i$ by 
Lemma~2.2 in \cite{Pedersen-Taylor(finiteness-obstruction)}. The claim now follows if we can prove that $j \circ \Phi \circ \chi$ is the 
identity. Let $[C_{\bullet}] \in \Sw^{\ch} ( \Gamma ; \IZ)$ be given. Without loss of generality we assume that
the complex is concentrated in non-negative degrees. There is an $m \geq 0$ such that the degrees of all non-vanishing homology groups lie 
in $\{ 0 ,1 , \dots , m \}$. We argue by induction over $m$. Let $m=0$ and put $M=H_0(C_{\bullet})$, i.e.\ $C_{\bullet}$
is a resolution of $M$.
Choose projective $\IZ \Gamma$-resolutions $P_{\bullet}^{\prime} \to TM$ and 
$P_{\bullet}^{\prime \prime} \to M/ T M$. Then $P_{\bullet}=P_{\bullet}^{\prime} \oplus
P_{\bullet}^{\prime \prime}$ is a resolution of $M$. Standard arguments produce chain maps 
$P_{\bullet}^{\prime} \to F_{\bullet} (TM)$, $P_{\bullet}^{\prime \prime} \to M/TM$ and $P_{\bullet} \to C_{\bullet}$ which are homology
isomorphisms. We have
\begin{eqnarray*}
& j \circ \Phi \circ \chi ( [C_{\bullet} ] ) = j \circ \Phi ( [M] ) = j ( [F_0] - [F_1] + [M/TM] ) = &\\
& [F_{\bullet}(TM)] + [M/TM]
= [P_{\bullet}^{\prime}] + [ P_{\bullet}^{\prime \prime} ]=[ P_{\bullet} ] = [C_{\bullet} ]. &
\end{eqnarray*}
Now let $m \geq 1$. Choose a projective $\IZ \Gamma$-resolution of $L_{\bullet} \to H_m( C_{\bullet} )$ and construct a $H_m$-isomorphism
$f:L_{\bullet} \to C_{\bullet}$. The map $f$ factorizes over its mapping cylinder and we have a short exact sequence
\[
0 \to L_{\bullet} \to \cyl (f) \to \cone (f) \to 0.
\]
The claim holds by the induction hypothesis for 
$[\cone (f) ]$ and by the argument for $m=0$ for $[L_{\bullet}]$ hence also for $[C_{\bullet} ]=[ \cyl (f) ]$.
\end{proof}

Let $\calc( \pt ; R \Gamma )$ denote the category of finitely generated free $R \Gamma$-modules.
Each $\IZ \Gamma$-module $M$ which is finitely generated free as a $\IZ$-module yields an exact functor 
$- \otimes_{\IZ} M : \calc( \pt ;  R \Gamma ) \to \calc ( \pt ; R \Gamma )$.
To check that this is well defined one
uses that for a free $R \Gamma$ module $F$ there is a non canonical isomorphism of $R \Gamma $-modules
between $F \otimes_{\IZ} M$ with the diagonal respectively with the left $\Gamma$-action.
Using Proposition~1.3.1 and 1.3.2~(4) in \cite{Waldhausen(1126)}  it is straightforward to check that the construction
leads to maps
\begin{eqnarray} \label{module-structure}
\Sw^{\fr} ( \Gamma ; \IZ ) \otimes_{\IZ} K_n ( R \Gamma ) \to K_n ( R \Gamma )
\end{eqnarray}
for $n \geq 0$. Replacing the one-point space $\pt$ by $\IR^n$ (compare Remark~\ref{remark-only-connected})
one obtains the corresponding construction for all $n \in \IZ$.
We use the isomorphism $i: \Sw^{\fr} ( \Gamma ; \IZ ) \to \Sw ( \Gamma ; \IZ )$ to define maps
\begin{eqnarray*} 
\Sw ( \Gamma ; \IZ ) \otimes_{\IZ } K_n ( R \Gamma ) \to K_n ( R \Gamma ).
\end{eqnarray*}
\begin{remark} The tensor product over $\IZ$ yields a ring structure on $\Sw^{\fr} ( \Gamma ; \IZ )$ and hence on
$\Sw ( \Gamma ; \IZ)$ and $\Sw^{\ch} ( \Gamma ; \IZ )$. The 
map (\ref{module-structure}) gives $K_n( R \Gamma )$ the structure of an $\Sw^{\fr} ( \Gamma ; \IZ )$-module.
\end{remark}

Let $\ch_{\hf} \calc ( \pt ; R \Gamma )$ denote the category of those  bounded below chain complexes which are
$\IZ \Gamma$-homotopy equivalent to a bounded complex of finitely generated free $\IZ \Gamma$-modules, compare
Subsection~\ref{appendix-homotopy-finite}.
Let $C_{\bullet}$ be a complex which represents an element in $\Sw^{\ch} ( \Gamma ; \IZ )$. 
\begin{lemma}
The functor $- \otimes_{\IZ} C_{\bullet}: \calc ( \pt ; R \Gamma ) \to \ch_{\hf} \calc ( \pt ; R \Gamma )$ is well defined.
\end{lemma}
\begin{proof}
Let $Q$ be a finitely generated free $\IZ \Gamma$-module. 
Since each $C_n$ is free as a $\IZ$-module there is a  non-canonical isomorphism between $Q \otimes_{\IZ} C_n$
equipped with the diagonal action and the same module with the right $\Gamma$-action. We see that $Q \otimes_{\IZ} C_{\bullet}$
is a complex of free $\IZ \Gamma$-modules. The crucial point is now to verify that this complex is homotopy equivalent to 
a bounded complex of finitely generated free $\IZ \Gamma$-modules. 
We argue by induction over the length of the interval in which the homology of $C_{\bullet}$ is concentrated.
We did this already in the proof of Proposition~\ref{all-three-isos} and will use the notation established there. 
Let $m=0$ and let $C_{\bullet}$ be a $\IZ \Gamma$-resolution of
the module $M=H_0 ( C_{\bullet} )$. 
Since $Q$ is free as a $\IZ$-module the sequence $0 \to Q \otimes_{\IZ} TM \to Q \otimes_{\IZ} M \to Q \otimes_{\IZ} M / TM \to 0$
is exact. 
The diagonal-action-versus-right-action argument used above shows that 
$Q \otimes_{\IZ} F_{\bullet}(TM)$ and $Q \otimes_{\IZ} M/TM$ are finitely generated free $\IZ \Gamma$-complexes. 
Using the horseshoe-lemma we obtain a resolution 
$P_{\bullet}=(Q \otimes_{\IZ } F_1 \to Q \otimes_{\IZ} F_0 \oplus Q \otimes_{\IZ} M/TM)$ of $M$ by a complex of finitely
generated free $\IZ \Gamma$-modules. Using standard arguments we can construct a $\IZ \Gamma$-chain map 
$P_{\bullet} \to Q \otimes_{\IZ} C_{\bullet}$ which is a homology isomorphism and hence a homotopy equivalence since both
complexes are complexes of projective $\IZ \Gamma$-modules. For the induction step one constructs 
a map $f:L_{\bullet} \to C_{\bullet}$
and a cylinder-cone sequence  as in the proof of  Proposition~\ref{all-three-isos}. One uses the induction hypothesis
and the horseshoe-lemma to construct a $\IZ \Gamma$-homotopy equivalence $P_{\bullet} \to C_{\bullet}$ with $P_{\bullet}$
bounded and degreewise finitely generated $\IZ \Gamma$-free.
\end{proof}

Since the inclusion $\calc( \pt ;  R \Gamma ) \to \ch_{\hf} \calc ( \pt ; R \Gamma )$ induces an isomorphism in $K$-theory 
(compare Lemma~\ref{lemma-inclusions-equivalences}) it is not difficult to check that we obtain maps
\[
\Sw^{\ch} ( \Gamma ; \IZ ) \otimes_{\IZ} K_n ( R \Gamma ) \to K_n ( R \Gamma ).
\]
Via the isomorphism $j$ from Proposition~\ref{all-three-isos} this action 
coincides with the action of $\Sw^{\fr} ( \Gamma ; \IZ )$ and $\Sw ( \Gamma ; \IZ )$. In particular 
$- \otimes_{\IZ} C_{\bullet}$
induces multiplication by $\sum (-1)^i [ H_i ( C_{\bullet} ) ]$.


\subsection{Some fibration sequences} \label{appendix-fibration-sequences}

Let $\cale$ be a morphism support condition on the $\Gamma$-space $X$
and let $\calf$, $\calf^{\prime}$, $\calf_0$ and $\calf_1$ be object support conditions.
Slightly abusing set-theoretical notation we define 
$\calf_0 \cap \calf_1=\{ F_0 \cap F_1 \; | \; F_0 \in \calf_0, \; F_1
\in \calf_1 \}$ and similarly $\calf_0 \cup \calf_1$.
Also we write $\calf_0 \subset \calf_1$ if for every $F_0 \in \calf_0$  
there is $F_1 \in \calf_1$ such that $F_0 \subset F_1$. If $\calf_0
\subset \calf_1$ and $\calf_1 \subset \calf_0$ we write $\calf_0 \simeq \calf_1$. We also use the corresponding 
notation for morphism support conditions.  

An object support condition $\calf$ is called {\em $\cale$-thickening closed} if
for every $F \in \calf$ and $E \in \cale$ there exists an $F^{\prime} \in \calf$ such that $F^E \subset F^{\prime}$.
Compare~\ref{subsubsec:thickenings} for notation. A typical example of such an $\calf$ is
$A^\cale = \{ A^E \; | \; E \in \cale \}$ for a subset $A \subset X$.
If $\calf$ and $\calf^{\prime}$ are $\cale$-thickening closed object support conditions,
then $\calc^{\Gamma} ( X, \cale , \calf \cap \calf^{\prime} )$ is 
a Karoubi filtration \cite{Karoubi(filtration)} of $\calc^{\Gamma} ( X , \cale , \calf)$ and we define 
$\calc^{\Gamma} ( X , \cale , \calf )^{> \calf^{\prime }}$
as the Karoubi quotient. For the definition of Karoubi filtrations and quotients we refer to \cite{Cardenas-Pedersen(karoubi)}.

\begin{lemma} \label{fibration-lemma} 
For the purpose of this lemma
we abbreviate 
\begin{eqnarray*}
(\calf ) & = & \calc^{\Gamma} ( X , \cale , \calf)  \\
(\calf )^{> \calf^{\prime}} & = & \calc^{\Gamma} ( X , \cale , \calf)^{> \calf^{\prime}}
\end{eqnarray*}
and we assume that all object support conditions which  occur are $\cale$-thickening closed.
\begin{enumerate}
\item \label{fibration-lemma-eins}
The sequence
$(\calf \cap \calf^{\prime}) \to (\calf) \to (\calf)^{> \calf^{\prime}}$
induces a fibration sequence in $K$-theory.
\item \label{fibration-lemma-zwei}
The square
\[
\xymatrix{
(\calf_0 \cap \calf_1 ) \ar[d] \ar[r] & (\calf_1) \ar[d] \\
(\calf_0 ) \ar[r] & (\calf_0 \cup \calf_1 )
         }
\]
induces a homotopy push-out square of spectra after applying $K$-theory.
\item \label{fibration-lemma-drei}
The sequence
$(\calf \cap \calf_0 )^{> \calf_1} \to (\calf)^{> \calf_1} \to (\calf)^{> \calf_0 \cup \calf_1}$
induces a fibration sequence in $K$-theory.
\end{enumerate}
\end{lemma}

\begin{proof}
The first statement is just the fact that a Karoubi filtration leads to a fibration sequence in $K$-theory,
compare \cite{Cardenas-Pedersen(karoubi)}. One can check that 
$(\calf_1)^{> \calf_0 \cap \calf_1} \to (\calf_0 \cup \calf_1 )^{> \calf_0}$
is an equivalence of categories. This yields \ref{fibration-lemma-zwei}.
The square 
\[
\xymatrix{
(\calf \cap \calf_0 \cap \calf_1 ) \ar[d] \ar[r] & (\calf \cap \calf_1) \ar[d] \\
(\calf \cap \calf_0 ) \ar[r] & (\calf \cap (\calf_0 \cup \calf_1) )
         }
\]
is a homotopy push-out square by \ref{fibration-lemma-zwei}. This square
maps to the homotopy push-out square whose vertical maps are identities and both whose 
horizontal maps are $(\calf \cap \calf_0) \to \calf$. The induced square of cofibers is again a homotopy push-out square and
its lower left hand corner is $(\calf \cap \calf_0)^{> \calf \cap \calf_0 }$ and hence contractible. Using 
\ref{fibration-lemma-eins} this yields the desired fibration sequence in \ref{fibration-lemma-drei}.
\end{proof}

\begin{example} \label{germs-at-infty}
Suppose that $\cale$ and $\calf$ are defined on $X \x [1,\infty)= X \times \IT$ and let
$\calf_{\IT} = \{ X \x [1,t_0] \; | \; t_0 \geq 1 \}$. Then we write $\calc^{\Gamma}(X \x \IT , \cale,\calf)^{\infty}$ for
$\calc^{\Gamma}(X \x \IT,\cale,\calf)^{> \calf_{\IT}}$ because we think of this 
category as being obtained by taking ``germs at infinity''. 
Let us  assume that the next remark applies to the inclusion $X = X \x \{ 1 \} \to X \x \IT$. 
Then the fibration sequence from Lemma~\ref{fibration-lemma}~\ref{fibration-lemma-eins}  
becomes
\[
\calc^{\Gamma} (X , \cale , \calf) \to \calc^{\Gamma} ( X \times \IT , \cale, \calf ) \to 
\calc^{\Gamma} ( X \times \IT , \cale, \calf )^{\infty}
\]
and we refer to it as the ``germs at infinity''-fibration.
\end{example}

\begin{remark} \label{proper-remark} 
Let $i_A : A \to X$ be the inclusion of a $\Gamma$-invariant subset.
Suppose that the morphism support condition $\cale$ satisfies the
following properness condition: for a compact $K \subset X$ and $E \in \cale$ the
closure of $K^E$ is again compact. Then 
\[
\calc^{\Gamma} ( A, i_A^{-1} \cale,i_A^{-1} \calf_0 )^{ > i_A^{-1} \calf}  
 \to \calc^{\Gamma} ( X , \cale , A^{\cale} \cap \calf_0 )^{> \calf} 
\]
is an equivalence of categories. 
(Here $\calf_0$ and $\calf$ are again assumed to be $\cale$-thickening closed object support conditions.)
Most morphism support conditions used in this paper are defined
on a locally compact metric space and contain a
global metric condition and are therefore proper in the above sense.
(In \cite{Bartels-Farrell-Jones-Reich(topology)} a weaker properness condition was used, but
this will not concern us here.)
\end{remark}

\begin{remark}[Mayer-Vietoris type results]
\label{Mayer-Viet}
Let $A$ and $B$ be $\Gamma$-invariant subsets of $X$ with $A \cup B = X$. 
Apply Lemma~\ref{fibration-lemma}~\ref{fibration-lemma-zwei} with $\calf_0= A^{\cale}$ and $\calf_1 = B^{\cale}$.
Suppose that $\cale$ is proper in the sense of
Remark~\ref{proper-remark} and that one can additionally show
\[
A^{\cale} \cap B^{\cale} \simeq (A \cap B)^{\cale}.
\]
Then one obtains a Mayer-Vietoris result, i.e.\ a homotopy push-out square involving the $K$-theories of the categories 
$\calc^{\Gamma}( Y, i^{-1} \cale)$, with $Y = A \cap B$, $A$, $B$ resp.\ $X$. 
(Here $i$ denotes in each case the relevant inclusion.) 
If $\calf'$ and $\calf$ are 
$\cale$-thickening closed object support conditions, then there is a similar 
homotopy push-out square for
the categories 
$\calc^{\Gamma}( Y, i^{-1} \cale, i^{-1} \calf)^{> i^{-1} \calf'}$. 
\end{remark}

\begin{lemma} \label{AFfequivalence}
Let $X$ be a space and $A \subset X$ a subspace. Let $\calf$ and $\calf_0$ be object support
conditions and $\cale^{\prime} \subset \cale^{\prime \prime}$ be morphism support conditions on $X$.
The map
\[
\calc(X , \cale^{\prime } , \calf_0 )^{> \calf \cup A^{\cale^{\prime}}} \to 
\calc(X , \cale^{\prime \prime} , \calf_0 )^{> \calf \cup A^{\cale^{\prime \prime}}} 
\]
induced by relaxing control from $\cale'$ to $\cale''$
is an equivalence of categories provided the following condition is satisfied.
\begin{quote}
For all $E^{\prime \prime} \in \cale^{\prime \prime}$ there exist 
$E^{\prime} \in \cale^{\prime}$ and $F \in \calf \cup A^{\cale'}$
such that 
\[
E^{\prime \prime} - E^{\prime} \subset F \times F.
\]
\end{quote}
\end{lemma} 

\begin{proof} 
Note that the condition implies $F \cup A^{\cale'} = F \cup A^{\cale''}$.  
But it also implies  that 
$\calc(X,\cale',\calf_0)^{> \calf \cup A^{\cale''}} \to 
 \calc(X,\cale'',\calf_0)^{> \calf \cup A^{\cale''}}$ is surjective on objects and 
bijective on morphism sets.
\end{proof}

The following concept was used to define strong control in Definition~\ref{definition-weak-strong-control}.

\begin{definition}[Control over a subset] \label{control-over}
Let $\cale$ be a control structure on the space $X$. Let $A \subset X$ be a subspace.
We define $\cale$-control over $A$ (a morphism-control condition on $X$) as follows.
A morphism $\phi$ is $\cale$-controlled over $A$ if there exists an $E \in \cale$
such that if $(x,y)$ lies in the support of $\phi$ and $x$ or $y \in A^{-E}$ then
$(x,y) \in E$.
\end{definition}

Here $A^{-E}=((A^c)^E)^c$ where $A^c$ denotes the complement of $A$ in $X$ and $A^E$ is the 
$E$-thickening of $A$ in $X$.


\typeout{----------------------- Literatur -------------------------}

\typeout{------------------------ THE END --------------------------}

\end{document}